\documentclass[11pt,leqno]{article}
\usepackage{graphicx}
\usepackage{amsfonts, amsthm}
\usepackage{amsxtra}
\usepackage[latin1]{inputenc}
\usepackage{fontenc}
\usepackage[dvips]{epsfig}
\usepackage[usenames]{color}
\begin{document}

\newtheorem{theorem}{Theorem}[section]
\newtheorem{definition}[theorem]{Definition}
\newtheorem{assumption}[theorem]{Assumption}
\newtheorem{example}[theorem]{Example}
\newtheorem{remark}[theorem]{Remark}
\newtheorem{lemma}[theorem]{Lemma}
\newtheorem{proposition}[theorem]{Proposition}
\newtheorem{corollary}[theorem]{Corollary}
\newtheorem{claim}[theorem]{Claim}

\numberwithin{equation}{section}
\numberwithin{theorem}{section}

\def\Z{\mathbb{Z}} 
\def\R{\mathcal{R}} 
\def\I{\mathcal{I}} 
\def\C{\mathbb{C}} 
\def\N{\mathbb{N}} 
\def\PP{\mathbb{P}} 
\def\Q{\mathbb{Q}} 
\def\L{\mathcal{L}} 
\def\ol{\overline} 
\def\bs{\backslash} 
\def\part{P} 

\begin{center}
{\bf  CONJUGACIES BETWEEN $P$- HOMEOMORPHISMS WITH
SEVERAL BREAKS}\footnote{ MSC: 37E10, 37C15, 37C40

Keywords and phrases: circle homeomorphism, break point, rotation
number, conjugation
} \\
\vspace{.25in} {\large {\sc Akhtam Dzhalilov\footnote{Turin
Politechnic University in Tashkent, Kichik halqa yuli, 17, 100095,
Tashkent,
 Uzbekistan. \quad \emph{E-mail}: a\_dzhalilov@yahoo.com},
Dieter Mayer\footnote{Institut f\"ur Theoretische Physik, TU
Clausthal, D-38678 Clausthal-Zellerfeld, Germany. \quad
\emph{E-mail}: dieter.mayer@tu-clausthal.de},

Utkir Safarov\footnote{ Institute of Mathematics,
National University of Uzbekistan, Tashkent, Uzbekistan.  \quad
\emph{E-mail}: safarovua@mail.ru}}}

\end{center}

\begin{abstract}

Let $f_{i},i=1,2$ be orientation preserving circle homeomorphisms
with a finite number of break points, at which the first derivatives
$Df_{i}$ have jumps, and  with identical irrational rotation number
$\rho=\rho_{f_{1}}=\rho_{f_{2}}.$ The jump ratio of $f_{i}$ at the
break point $b$ is denoted by $\sigma_{f_{i}}(b)$, i.e.
$\sigma_{f_{i}}(b):=\frac{Df_{i}(b-0)}{Df_{i}(b+0)}$. Denote by
$\sigma_{f_{i}}, i=1,2,$
 the total jump ratio given by the product over all break points $b$ of the jump ratios $\sigma_{f_{i}}(b)$ of $f_{i}$.
 We prove, that for circle homeomorphisms
$f_{i}, i=1,2$, which are $C^{2+\varepsilon}, \varepsilon>0$,  on each
interval of continuity of $Df_{i}$ and whose total jump ratios
$\sigma_{f_{1}}$ and $\sigma_{f_{2}}$ do not coincide, the
congugacy between  $f_{1}$ and  $f_{2}$ is a singular function.
\end{abstract}

\section{Introduction}

Let $f$ be an orientation preserving homeomorphism of the circle
$S^{1}\equiv\mathbb{R}^{1}/\mathbb{Z}^{1}$ with
lift $F:\mathbb{R}^{1}\rightarrow{\mathbb{R}^{1}}$, which is
continuous, strictly increasing and fulfills $F(x+1)=F(x)+1$,
$x\in\mathbb{R}$. The circle homeomorphism $f$ is then defined by
$f(x)=F(x)$ $(mod\,1)$, $x\in{S^{1}}.$

 Denjoy's classical  theorem \cite{De1932} states, that a
circle diffeomorphism $f$ with irrational rotation number
$\rho=\rho_{f}$ and such that ${\log}Df$ is of
bounded variation, is conjugate to the linear rotation $f_{\rho}$,
that is, there exists a homeomorphism $\varphi$ of the circle with
$f=\varphi^{-1} \circ f_{\rho} \circ \varphi.$

It is well known that a circle homeomorphisms $f$ with irrational
rotation number $\rho$ is strictly ergodic, i.e. it has a unique
$f$- invariant probability measure $\mu_{f}$. A remarkable fact is
then that the conjugacy $\varphi$ can be defined by
$\varphi(x)=\mu_{f}([0,x])$ (see \cite{CFS1982}), which shows, that
the regularity properties of this conjugacy $\varphi$ imply the
corresponding properties of the density of the absolutely continuous
invariant measure $\mu_{f}$. The problem of smoothness of the
conjugacy of smooth diffeomorphisms is by now very well understood
(see for instance
\cite{Ar1961},\cite{Mo1966},\cite{He1979},\cite{KO1989},\cite{KS1989},\cite{Yo1984}).
Notice, that for a sufficiently smooth circle diffeomorphism with a
typical irrational rotation number its invariant measure is
absolutely continuous with respect to Lebesque measure (see
\cite{KO1989},\cite{KS1989}).

A natural extension of circle diffeomorphisms are piecewise smooth
homeomorphisms with break points or shortly, the class of P-homeomorphisms.

 This class of \textbf{P-homeomorphisms}
 consists of orientation preserving circle homeomorphisms
$f$ which are differentiable except at a finite number of break points, at which the
one-sided positive derivatives  $Df_{-}$ and $Df_{+}$ exist, which do not coincide and for which there exist constants  $0<c_{1}<c_{2}< \infty$,
such that

$\bullet$\,\,$c_{1}<Df_{-}(x_{b})<c_{2}$ and
$c_{1}<Df_{+}(x_{b})<c_{2}$ for all $x_{b}\in {B(f)}$, the set
 of break points of $f$ in $S^{1}$;

$\bullet$\,\, $c_{1}<Df(x)<c_{2}$
 for all $x\in S^{1}\backslash{B(f)}$;

$\bullet$\,\,$\log Df$ has bounded variation in $S^{1}$ i.e.
$v:=var_{S^{1}}\log Df<\infty$.

The ratio $\sigma_{f}(x_b)=\frac{Df_{-}(x_{b})}{Df_{+}(x_{b})}$ is
called the \textbf{jump ratio} of $f$ at  $x_b$  or, for short, the
$f$-\textbf{jump}. The product of all jump ratios  is called the
\textbf{total jump} of $f$ and denoted by $\sigma_{f}.$ Notice, that
 Denjoy's result can be extended to $P$-homeomorphisms  with
irrational rotation numbers, its precise formulation will be given
later.

 The regularity properties of the invariant measures of
P-homeomorphisms are quite different from those of diffeomorphisms
(see \cite{DK1998},\cite{L2004}, \cite{DL2006}, \cite{DLM2009},
\cite{DMS2012}). Dzhalilov, Mayer and Safarov  proved in
\cite{DMS2012}, that the invariant measures of piecewise
$C^{2+\varepsilon}$ P-homeomorphisms $f$ with non trivial total jump
$\sigma_f$ and with irrational rotation number are singular w.r.t.
Lebesgue measure. In this case the conjugacy $\varphi$ between $f$
and the linear rotation $f_{\rho}$ is a singular function. Here then
arises naturally the problem of regularity of the conjugacy between
two circle maps with identical irrational rotation numbers and with
break points.
 This is the so called rigidity problem for circle
homeomorphisms with  break points.  The case of two circle maps with
one break point and  the same jump ratio were studied in detail by
K. Khanin and D. Khmelev \cite{KhKm2003}, K. Khanin and A. Teplinsky
\cite{KT2013}. To formulate their result, let $\rho=1/\left(k_{1}+
1/\left( k_{2}+...+1/\left(k_{n}+...\right) \right)
\right):=[k_1,k_2,\ldots,k_n,\ldots]$ be the continued fraction
expansion of the irrational rotation number $\rho.$

Define
$$
M_o=\{\rho:\exists C>0, \forall n\in \mathbb{N},\, k_{2n-1}\leq
C\},\,\,\,\, M_e=\{\rho:\exists C>0, \forall n\in \mathbb{N},\,
k_{2n }\leq C\}.
$$
Then  K. Khanin and A. Teplinskii proved in \cite{KT2013}
\begin{theorem}\label{Tep Kh}
Let  $f_{i}\in C^{2+\alpha}(S^{1}\backslash\{{b_{i}}\}),  \ i=1,2$,
$\alpha>0$, be  circle homeomorphisms with one break point, the same
jump ratio $\sigma$ and the same irrational rotation
number $\rho\in (0, 1).$ If either $\sigma> 1$ and $\rho\in M_e$ or
$\sigma<1$ and $\rho\in M_o$, then the map $h$ conjugating the
homeomorphisms $f_ 1$ and $f_2$ is a $C^1$-diffeomorphism.
\end{theorem}
In the case of homeomorphisms with different jump ratios the
following theorem was proved by A. Dzhalilov, H. Akin and S. Temir in
\cite{DHT2010}:
\begin{theorem}\label{Dz.H.T.}
Let  $f_{i}\in C^{2+\alpha}(S^{1}\backslash\{{b_{i}}\}),  \ i=1,2$,
$\alpha>0$,  be  circle homeomorphisms with one break point and
 different jump ratio but the same irrational rotation number
$\rho\in (0, 1).$ Then the map $h$ conjugating the homeomorphisms
$f_ 1$ and $f_2$ is a singular function.
\end{theorem}
Now consider two piecewise-smooth circle homeomorphisms  $f_{1}$ and
$f_{2}$  with $m$ $(m\geq2)$ break points and the same irrational
rotation number. Denote by $B(f_{1})$ and $ B(f_{2})$ the set of
break points of $f_{1}$ and  $f_{2}$ respectively.
 \begin{definition}\label{defKO}
The homeomorphisms $f_{1},$ $f_{2}$ are said to be \textbf{break
point equivalent} if there exists a topological conjugacy $\psi_{0}$
such that
\begin{enumerate}
  \item [(1)] $\psi_{0}(B(f_{1}))=B(f_{2})$;
  \item [(2)] $\sigma_{f_{2}}(\psi_{0}(b))=\sigma_{f_{1}}(b),$ for all $b\in B(f_{1}).$
\end{enumerate}
\end{definition}
The rigidity problem for break point equivalent
$C^{2+\alpha}$-homeomorphisms $f$ with trivial total jumps $\sigma_f =1$
 was studied by K. Cunha and D. Smania in
\cite{CS2012}. It was proven there that any two such homeomorphisms
 fulfilling certain combinatorial conditions are $C^{1}$-conjugated.
  The main idea
of their proof is to consider piecewise-smooth circle homeomorphisms
as
 generalized interval exchange transformations.
The case of non break point equivalent homeomorphisms with two break
points was studied by H. Akhadkulov, A. Dzhalilov and D. Mayer in
\cite{ADM2012}.
 Their main result  is the following theorem:
\begin{theorem}\label{ADM }  Let  $f_{i}\in C^{2+\alpha}(S^{1}\backslash\{{a_{i},b_{i}}\}), i=1,2 $
be circle homeomorphisms with two break points $ a_{i},b_{i}.$
Assume that
\begin{enumerate}
  \item [(1)] their rotation numbers $\rho_{f_{i}},$ $i=1,2$ are irrational and coincide
  i.e. $\rho_{f_{1}}=\rho_{f_{2}}=\rho$;
  \item [(2)] there exists a bijection  $\psi$ such that $\psi(B(f_{1}))=B(f_{2})$;
  \item [(3)] $\sigma_{f_{1}}=\sigma_{f_{1}}(a_{1})\sigma_{f_{1}}(b_{1})\neq\sigma_{f_{2}}
  =\sigma_{f_{2}}(a_{2})\sigma_{f_{2}}(b_{2})$.
  \end{enumerate}
Then the map $h$ conjugating $f_1$ and $f_2$ is a singular function.
\end{theorem}

  In the present paper we study the conjugacy $h$
of two piecewise smooth circle homeomorphisms $f_{1}$ and
$f_{2}$ with an arbitrary finite number of  break points .

Our main result is the following theorem:
\begin{theorem}\label{DMS}
Let $f_{i},$  $i=1,2,$ be P-homeomorphisms with the same irrational
rotation number  $\rho=\rho_{f_{1}}=\rho_{f_{2}}$.
   Assume, that
\begin{enumerate}
\item [(1)] $f_{i}$, $i=1,2$ is
$C^{2+\alpha},\,\alpha>0,$ on each interval of continuity  of
$Df_i$;
 \item [(2)]\,the total jumps of $f_{1}$ and $f_{2}$ do not
 coincide i.e.

 $\sigma_{f_{1}}=\prod\limits_{b\in B(f_1)}\sigma_{f_{1}}(b)\neq
  \sigma_{f_{2}}=\prod\limits_{b\in B(f_2)} \sigma_{f_{2}}(b).$
\end{enumerate}
Then the  map  $h$ conjugating $f_{1}$ and $f_{2}$ is a singular
function.
\end{theorem}

\section{Preliminaries and Notations}

Let $f$ be an orientation-preserving circle homeomorphism  with lift
$F$. The important characteristic of homeomorphism $f$ is the
rotation number  defined by
  $$\rho_{f} :=
\underset{{n\rightarrow\infty}}{\lim}
\frac{F^{n}(x)}{n}\,\,\,(mod1).$$ Here and below, $F^{n}$ denotes
the $n$-th iteration of $F$. Suppose the rotation number $\rho_{f}$
is irrational.Then it can be uniquely represented as a continued
fraction i.e. $\rho_f: = [k_{1}, k_{2},..., k_{n},... ]$.
 Define
$\frac{p_n}{q_n}: = [k_{1}, k_{2},..., k_{n}], n\geq 1$ the
convergent of $\rho_{f}$. Their denominators   $q_n$ satisfy the
recurrence relation: $q_{n+1} = k_{n+1}q_{n}+q_{n-1}, n\geq 1$, with
the initial conditions $q_{0}=1$ and $q_{1}=k_{1}$.

Fix a point $x_0\in S^{1}$. Its positive orbit
$\{x_i=f^{i}(x_{0}),\,\,i=0,1,2...\}$  defines a sequence of
natural partitions of the circle: denote by $\Delta_{0}^{(n)}(x_{0})$
 the closed interval in $S^{1}$ with endpoints $x_{0}$ and
$x_{q_{n}}=f^{q_{n}}(x_{0}).$ Notice, that for $n$ odd the point
$x_{q_{n}}$ is to the left of $x_0$, and for $n$ even it is to its
right. Denote by
$\Delta_{i}^{(n)}(x_{0})=f^{i}(\Delta_{0}^{(n)}(x_{0})), i\geq 1$,
the iterates of the interval $\Delta_{0}^{(n)}(x_{0})$ under $f$. It
is well known, that the set $\eta_{n}(x_{0})$ of
intervals with mutually disjoint interiors defined as
$$
\eta_{n}(x_0)=\{\Delta_{i}^{(n-1)}(x_{0}),\,\,0\leq i\leq
q_{n}-1\}\cup\{\Delta_{j}^{(n)}(x_{0}),\,\,0\leq j\leq q_{n-1}-1\}
$$
determines a partition of the circle for any $n$. The partition
$\eta_{n}(x_{0})$ is called the $n$-th
\textbf{dynamical partition} of the point $x_{0}$. Proceeding from
$\eta_{n}(x_{0})$ to $\eta_{n+1}(x_{0})$ all the
intervals $\Delta_{j}^{(n)}(x_{0}),\,\,0\leq j\leq
q_{n-1}-1$, are preserved, whereas each of the intervals
$\Delta_{i}^{(n-1)}(x_{0})$ is partitioned into $k_{n}+1$
subintervals belonging to $\eta_{n+1}(x_{0})$ such
that
$$ \Delta_{i}^{(n-1)}(x_{0})=
\Delta_{i}^{(n+1)}(x_{0})\cup\bigcup\limits_{s=0}^{k_{n+1}-1}\Delta_{i+q_{n-1}+sq_{n}}^{(n)}(x_{0}).
$$
Obviously one has $\eta_{1}(x_{0})\leq
\eta_{2}(x_{0})\leq ...\leq \eta_{n}(x_{0})\leq ...$.
\begin{definition}\label{COM} Let $K>1$ be a constant. We call two intervals
$I_{1}$ and $I_{2}$ of $S^{1}$ $K$-comparable, if the
inequality  $ K^{-1}\ell(I_{2})\le \ell(I_{1})\le K
\ell(I_{2})$ holds.
\end{definition}

Following \cite{KO1989} we recall
definition.
\begin{definition}\label{SMALL} An interval $I=[\tau,t]\subset S^{1}$ is said to
be $q_{n}$-small, and its endpoints $q_{n}$-close, if the intervals
$f^{i}(I),\; 0\le i\leq q_{n}-1$, are pairwise disjoint (except for
endpoints).
\end{definition}
It follows from the structure of the dynamical partitions that an
interval $I=[\tau,t]$ is $q_{n}$-small if and only if either
$\tau\prec t \preceq f^{q_{n-1}}(\tau)$, or $f^{q_{n-1}}(t)\preceq
\tau\prec t$.
\begin{lemma}\label{FIN} Let $f$ be a class $P$-homeomorphism with a finite number of break points and irrational
rotation number $\rho=\rho_{f}$. If the interval $I=(x,y)\subset
S^{1}$ is $q_{n}$- small and $f^{s}(x),\,f^{s}(y)
\not\in B(f)$ for all $0\leq s<q_{n}$, then
for any $k\in [0,q_{n}]$ the following inequality holds:
\begin{equation}\label{FLK}
e^{-v}\leq \frac{ Df^{k}(x)}{Df^{k}(y)}\leq e^{v},
\end{equation}
where $v$ is total variation of $\log Df$ in $S^{1}$.
\end{lemma}

\textbf{Proof of Lemma \ref{FIN}.} Take any two $q_{n}$-close points
$x,y\in S^{1}$ and $0\leq k\leq q_{n}-1$. Denote by $I$ the open
interval with endpoints $x$ and $y$. Because the intervals
$f^{i}(I),\,\,0\leq i\leq q_{n}-1$ are disjoint, we obtain
$$
|\ln Df^{k}(x)-\ln Df^{k}(y)|\le \sum\limits_{j=0}^{k-1}|\ln
Df(f^{j}(x))-\ln Df(f^{j}(y))|\le v,
$$
from which inequality (\ref{FLK}) follows immediately.

Using  Lemma \ref{FIN} the
following lemma can be proven which plays a key role in the  study of the metrical
properties of homeomorphisms.
\begin{lemma}\label{DEN} Suppose the circle homeomorphism $f$ satisfies the conditions of Lemma \ref{FIN}. Then for any
$y_{0}$ with  $y_s:=f^{s}(y_{0})\not\in B(f)$,\,\, for all $  s\in[0,q_{n})$ the following inequality holds:
\begin{equation}\label{PLF}
e^{-v}\le \prod\limits_{s=0}^{q_{n}-1}Df(y_{s})\le e^{v}.
\end{equation}
\end{lemma}

Inequality (\ref{PLF}) is called the Denjoy's inequality. It follows
from Lemma \ref{DEN} that the intervals of the dynamical partition
$\eta_{n}(x_{0})$ have exponentially small lengths.
Indeed one finds
\begin{corollary}\label{LDB} Let $\Delta^{(n)}$ be an arbitrary element of
the dynamical partition $\eta_{n}(x_0)$. Then
\begin{equation}\label{EDL}
 \ell(\Delta^{(n)})\leq const \lambda^n
\end{equation}
where $\lambda=(1+e^{-v})^{-\frac{1}{2}}<1$.
\end{corollary}
\begin{definition}\label{THTT} Two homeomorphisms $f_{1}$ and $f_{2}$
of the circle are said to be topologically equivalent if there
exists a homeomorphism $\varphi:S^{1}\to S^{1}$ such that
$\varphi(f_{1}(x)) = f_{2}(\varphi(x))$ for any $x\in S^{1}$.
\end{definition}
 The homeomorphism $\varphi$ is called a \textbf{conjugacy}.
  Corollary
\ref{LDB}  implies the following generalization of the classical
Denjoy theorem:
 \begin{theorem}\label{Denjoy} Suppose that a
homeomorphism $f$ satisfies the conditions of Lemma \ref{FIN}. Then
the homeomorphism $f$ is topologically conjugate to the linear
rotation $f_{\rho}$.
\end{theorem}

 In the proof of our main theorem the
tool of {\bf{cross- ratio}} plays a key role.
\begin{definition}\label{Cross} The {\bf cross-ratio} of four numbers
$(z_{1},z_{2},z_{3},z_{4}),\; z_{1}<z_{2}<z_{3}<z_{4}$, is the
number
$$
Cr(z_{1},z_{2},z_{3},z_{4})=\frac{(z_{2}-z_{1})(z_{4}-z_{3})}{(z_{3}-z_{1})(z_{4}-z_{2})}.
$$
\end{definition}
\begin{definition}\label{Dist} Given four real numbers $(z_{1},z_{2},z_{3},z_{4}) $ with
$ z_1<z_2<z_3<z_4$ and a strictly increasing function $F:\; \mathbb{R}^1\to
\mathbb{R}^1$. The distortion of their cross-ratio under $F$ is given by
$$
Dist(z_{1},z_{2},z_{3},z_{4};F)=\frac{Cr(F(z_1),F(z_2),F(z_3),F(z_4))}{Cr(z_{1},z_{2},z_{3},z_{4})}.
$$
\end{definition}

For $m\geq 3$ and $z_{i}\in S^{1},\,1\leq i\leq m$,  suppose that
$z_{1}\prec{z_{2}}\prec...\prec{z_{m}}\prec{z_{1}}$ (in the sense of
the ordering on the circle). Then we set $\hat{z}_{1}:=z_{1}$ and

\[\hat{z}_{i}:=\left\{\begin{array}{ll}z_{i}, & \mbox{if  $z_{1}\prec z_{i}\prec 1$},\\1+z_{i}, &
 \mbox{if
 $0\prec{z_{i}}\prec z_{1}$}.\end{array}\right.\] for $2\leq i \leq
 m$.

Obviously, $\hat{z}_{1}<{\hat{z}_{2}}<{...}<{\hat{z}_{m}}$. The
vector $(\hat{z}_{1},\hat{z}_{2},...,\hat{z}_{m})$ is called the
lifted vector of $(z_{1},z_{2},...,z_{m})\in (S^{1})^{m}$.

Let $f$ be a circle homeomorphism with lift $F$. We define the
cross-ratio distortion of
$(z_{1},z_{2},z_{3},z_{4}),\,z_1\prec{z_{2}}\prec{z_{3}}\prec{z_{4}}\prec{z_{1}}$
with respect to $f$ by
$Dist(z_{1},z_{2},z_{3},z_{4};f)=Dist(\hat{z}_{1},\hat{z}_{2},\hat{z}_{3},\hat{z}_{4};F)$,
where $(\hat{z}_{1},\hat{z}_{2},\hat{z}_{3},\hat{z}_{4})$ is the
lifted vector of $(z_{1},z_{2},z_{3},z_{4})$. We need the following
\begin{lemma}\label{DAKX1}(see \cite{DK1998}). Let $z_{i}\in S^{1},\,i=1,2,3,4,\,
z_1\prec{z_{2}}\prec{z_{3}}\prec{z_{4}}.$ Consider a circle
homeomorphism $f$ with $f\in
C^{2+\varepsilon}([z_{1},z_{4}]),\,\varepsilon>0,$ and $Df(x)\geq
const>0$ for $x\in [z_{1},z_{4}].$ Then there is a positive constant
$C_{1}=C_{1}(f)$ such that
$$
\mid Dist(z_{1},z_{2},z_{3},z_{4};f)-1\mid\leq
C_{1}|\hat{z}_{4}-\hat{z}_{1}|^{1+\varepsilon},
$$
where $(\hat{z}_{1},\hat{z}_{2},\hat{z}_{3},\hat{z}_{4})$ is the
lifted vector of $(z_{1},z_{2},z_{3},z_{4})$.
\end{lemma}

We next consider the case where the interval $[z_{1},z_{4}]$ contains a
break point $b$ of the homeomorphism $f$. More precisely, suppose
 $b\in [z_{1},z_{2}]$. Let $\sigma_{f}(b)$ be the jump of $f$ at
$b$. We define numbers $\alpha,\,\beta,\,\tau,\,\xi$ and $z$ as
follows:
$$\alpha:=\hat{z}_{2}-\hat{z}_{1},\,\,
\beta:=\hat{z}_{3}-\hat{z}_{2}, \,\,\tau:=\hat{z}_{2}-\hat{b},
\,\,\xi:=\frac{\beta}{\alpha},\,\,z:=\frac{\tau}{\alpha}.
$$
where $(\hat{z}_{1},\hat{b},\hat{z}_{2},\hat{z}_{3})$ is the lifted
vector of $(z_1,b,z_2,z_3)$.

In what follows we shall need the following lemma.

\begin{lemma}\label{DAKX2}(see \cite{DK1998}). For the circle  homeomorphism $f$ with $f\in
C^{2}([z_{1},z_{4}]\backslash \{b\}),$ and  $Df(x)\geq const>0$ for
$x\in [z_{1},z_{4}]\backslash \{b\}$ one has $$
\left| Dist(z_{1},z_{2},z_{3},z_{4};f)-
\frac{[\sigma_{f}(b)+(1-\sigma_{f}(b))z](1+\xi)}{\sigma_{f}(b)+(1-\sigma_{f}(b))z+\xi}\right|\leq
C_{2}|\hat{z}_{4}-\hat{z}_{1}|,
$$
where the constant $C_{2}>0$ depends only on $f$.
\end{lemma}

\section{On $q_{n}$-preimages of break points}

Let  $f_{1}$ and $f_{2}$ be $P$-homeomorphisms with identical
irrational rotation numbers $\rho=\rho_{f_{1}}=\rho_{f_{2}}$. Denote
by
 $B(f_{1})=\{b_{1}^{(i)},1\leq i\leq m_{1}\}$ and
$B(f_{2})=\{b_{2}^{(j)}, 1 \leq j\leq m_{2}\}$ the sets of all break
points of $f_{1}$ and $f_{2}$, respectively. Take two copies of the
circle on which the homeomorphisms $f_{1}$ and $f_{2}$ act
respectively. Denote by $\varphi_{i},\,i=1,2$  the conjugacies
between $f_{i}$ and $f_{\rho}$ i. e. $\varphi_{1}\circ
f_{1}=f_{\rho}\circ \varphi_{1}$ and $\varphi_{2}\circ
f_{2}=f_{\rho}\circ \varphi_{2}.$ It is easy to check that the
homeomorphisms $f_{1}$ and $f_{2}$ are then conjugated by
$h=\varphi_{2}\circ \varphi_{1}^{-1}$ i. e. $h\circ
f_{1}(x)=f_{2}\circ h(x), \forall x\in S^{1}.$   For $x_{0}\in
S^{1}$ let $\eta_{n}(x_0)$ be its $n$-th dynamical
partition.  Put $t_{0}:=h(x_{0})$ and  consider the
dynamical partition $\tau_{n}(t_0)$ of $t_{0}$ on the second circle
determined by the homeomorphism $f_{2}$ i.e. $$
\tau_{n}(t_0)=\{C_{i}^{(n-1)}(t_0),\,\,0\leq i\leq
q_{n}-1\}\cup\{C_{j}^{(n)}(t_0),\,\,0\leq j\leq q_{n-1}-1\}.
$$
with $C_0^{(n)}(t_0)$ the closed interval with endpoints $t_0$ and
$f_2^{q_n}(t_0)$. Chose an odd natural number  $n=n(f_{1},f_{2})$
such that the $n$-th renormalization neighborhoods
$[x_{q_{n}},x_{q_{n-1}}]$ and $[t_{q_{n}},t_{q_{n-1}}]$ do not
contain any break point of $f_{1}$ and $f_{2}$ respectively. Since
the identical rotation number $\rho$ of $f_{1}$ and $f_{2}$ is
irrational, the order of the points on the orbit
$\{f_{1}^k(x_0),\,\,k\in \mathbb{Z}\}$ on the first circle will be
precisely the same as the one for the orbit $\{f_{2}^k(t_0),\,\,k\in
\mathbb{Z}\}$ on the second one. This together with the relation
$h(f_{1}(x))=f_{2}(h(x))$ for $x\in S^1$ implies that
\begin{equation}\label{HDC}
h(\Delta_{i}^{(n-1)})=C_{i}^{(n-1)},\,\,0\leq i\leq
q_{n}-1,\,\, h(\Delta_{i}^{(n)})=C_{j}^{(n)},\,\,0\leq j\leq
q_{n-1}-1.
\end{equation}

The structure of the dynamical partitions implies that
$\overline{b}_{1}^{(i)}(n)=f_{1}^{-l_{1}^{(i)}}(b_1^{(i)})\in[x_{q_n},x_{q_{n-1}}],\,\,1\leq
i \leq m_1$, where $l_{1}^{(i)}\in (0,q_{n-1})$ if
$\overline{b}_{1}^{(i)}(n)\in [x_{q_n},x_{0}]$, and  $l_{1}^{(i)}\in
(0,q_{n})$ if $\overline{b}_{1}^{(i)}(n)\in [x_{0},x_{q_{n-1}}]$.
Also $\overline{b}_{2}^{(j)}(n)=f_{2}^{-l_{2}^{(j)}}(b_2^{(j)})\in
[t_{q_n},t_{q_{n-1}}],\,\,1\leq j\leq m_2$ where $l_{2}^{(j)}\in
(0,q_{n-1})$ if $\overline{b}_{2}^{(j)}(n)\in [t_{q_n},t_{0}]$, and
$l_{2}^{(j)}\in (0,q_{n})$ if $\overline{b}_{2}^{(j)}(n)\in
[t_{0},t_{q_{n-1}}]$. The points $\overline{b}_{1}^{(i)}(n)$ and
$\overline{b}_{2}^{(j)}(n)$ are called the $q_n$-\textbf{preimages}
of the break points $b_{1}^{(i)}$ and $b_{2}^{(j)}$. Denote by
$\overline{B}^{(n)}(f_{i}), i=1,2$,  the sets of
$q_n$-preimages in the renormalization intervals
$[x_{q_n},x_{q_{n-1}}]$ and $[t_{q_n},t_{q_{n-1}}]$ of the sets
$B(f_{1})$ and $B(f_{2})$, respectively. Then the number of points
in $\overline{B}^{(n)}(f_{i}), i=1,2$, is not
greater than $m_{i}$. Next, consider the set
 $h^{-1}(\overline{B}^{(n)}(f_{2}))$. Using the relations (\ref{HDC}) we
find that $h^{-1}(\overline{B}^{(n)}(f_{2}))\subset
[x_{q_n},x_{q_{n-1}}]$. Notice, that the number of elements of the
set $\overline{B}^{(n)}(f_{1})\cup
h^{-1}(\overline{B}^{(n)}(f_{2}))$ is bounded by $m_1+m_2$.

We set \begin{equation}\label{BMM}
B_{m_{1},m_{2}}^{(n)}=\{x_{q_{n}},\,x_{0},\,x_{q_{n-1}}\}\cup\overline{B}^{(n)}(f_{1})\cup
h^{-1}(\overline{B}^{(n)}(f_{2})),\,\,
d_{n}=\ell([x_{q_{n}},x_{q_{n-1}}]).
\end{equation}
Let $m_{0}\in N,\,m_{0}>m_{1}+m_{2}+3.$ For every $l\geq 0$ we
define a partition $D_{l}^{(n)}$ of the interval
$[x_{q_{n}},x_{q_{n-1}}]$ using the points
$t_{s}=x_{q_{n}}+m_{0}^{-(l+1)}{d_{n}}s ,\,s=0,1,...,m_{0}^{l+1}$.
Obviously, the length of every such interval $I_{l}^{(n)}$ of
$D_{l}^{(n)}$ is equal to $m_{0}^{-(l+1)}{d_{n}}$. When passing from
$D_{l}^{(n)}$ to $D_{l+1}^{(n)}$, every interval of $D_{l}^{(n)}$ is
divided into $m_0$ intervals $D_{l+1}^{(n)}$.
\begin{definition} An
interval of $D_{l}^{(n)}$ that does not contain any elements of
$B_{m_{1},m_{2}}^{(n)}$ is called an $l^{(n)}$-empty interval.
Otherwise it is called an $l^{(n)}$-occupied interval.
\end{definition}
Since the number of intervals of $D_{0}^{(n)}$ is greater than the
number of elements of $B_{m_{1},m_{2}}^{(n)}$, there exists at least
one $0^{(n)}$-empty interval. Furthermore, every $l^{(n)}$-occupied
interval of $D_{l}^{(n)}$ contains at least one $(l+1)^{(n)}$-empty
interval of $D_{l+1}^{(n)}$. Note that the leftmost and rightmost
intervals of $D_{l}^{(n)}$ contain $x_{q_{n}}$ and $x_{q_{n+1}}$,
respectively. This means that these extreme intervals are
$l^{(n)}$-occupied for any $l\geq 0$. Removing all  $l^{(n)}$-empty
intervals from the interval $[x_{q_{n}}, x_{q_{n+1}}]$, we obtain a
natural partition of $B_{m_{1},m_{2}}^{(n)}$ into non-empty disjoint
parts. Denote this partition by $\Gamma_{l}^{(n)}$. Between two
elements of the partition $\Gamma_{l}^{(n)}$ lies at least one
$l^{(n)}$-empty interval. Removing all $l^{(n)}$-occupied intervals
from the interval $[x_{q_{n}}, x_{q_{n+1}}]$ we obtain the set
$V_{l}^{(n)}$  of intervals.

The structure of the set
$\overline{B}^{(n)}(f_{1})\cup
h^{-1}(\overline{B}^{(n)}(f_{2}))$ in $[x_{q_{n}}, x_{q_{n+1}}]$ is
given by
\begin{theorem}\label{LG} Let  $f_{1}$,
 $f_{2}$ be $P$-homeomorphisms with a finite number of break points
  and with identical irrational rotation numbers. Suppose that their total
  jumps do not coincide.
 For any
positive integer $r$ there exists a number $s_{0}=s_{0}(r,n),\,
0\leq{s_{0}}\leq{r(m_1+m_2+1)}$, such that
\newline
 1) $\max\limits_{x,y\in E_{s_{0}}^{(n)},x\prec y}\ell([x,y])\leq
2m_{0}^{-(s_{0}+r)}d_{n}$ for every $E_{s_{0}}^{(n)}\in
\Gamma_{s_{0}}^{(n)};$\newline
 2) $\ell(I)\geq m_{0}^{-s_{0}}d_{n},$
for all $I\in V_{s_{0}}^{(n)};$\newline
 3) there exists at least one
element $\widetilde{E}_{s_{0}}^{(n)}$ of the partition
$\Gamma_{s_0}^{(n)}$ such that
\begin{equation}\label{PROD}\prod\limits_{b_{1}:\overline{b}_{1}(n)\in
\widetilde{E}_{s_{0}}^{(n)}}\sigma_{f_{1}}(b_{1})\neq
\prod\limits_{b_{2}:\overline{b}_{2}(n)\in
h(\widetilde{E}_{s_{0}}^{(n)})}\sigma_{f_{2}}(b_{2}).
\end{equation}
\end{theorem}

\textbf{Proof of Theorem \ref{LG}.} Consider the partitions
$D_{0}^{(n)},\,D_{r}^{(n)},\,...,\,D_{r(m_1+m_2+1)}^{(n)}$ and the
partitions $\Gamma_{0}^{(n)},\,\Gamma_{r}^{(n)},...,
\Gamma_{r(m_1+m_2+1)}^{(n)}$ of the set $B_{m_{1},m_{2}}^{(n)}$
generated by them. Let $|\Gamma_{l}^{(n)}|$ denote the number of
elements of the partition $\Gamma_{l}^{(n)}.$ For proving the first
two assertions of the theorem it is sufficient to show that
$|\Gamma_{s_{0}}^{(n)}|=|\Gamma_{s_{0}+r}^{(n)}|,$ for some $s_{0}.$

 It follows from the structure of
the partitions $\Gamma_{l}^{(n)}$ that $|\Gamma_{l}^{(n)}|\leq
|\Gamma_{l+1}^{(n)}|$ for any $l\geq 0$. In particular,
$|\Gamma_{0}^{(n)}|\leq |\Gamma_{r}^{(n)}|\leq ...\leq
|\Gamma_{r(m_1+m_2+1)}^{(n)}|$. Then two
cases are possible: either
$|\Gamma_{rt_{0}}^{(n)}|=|\Gamma_{r(t_{0}+1)}^{(n)}|$ for some
$t_{0}=t_{0}(n),\,0\leq t_{0} \leq m_1+m_2,$ or $|\Gamma_{0}^{(n)}|<
|\Gamma_{r}^{(n)}|< ...< |\Gamma_{r(m_1+m_2+1)}^{(n)}|$. In the
first case we set $s_{0}=rt_{0}$. If $|\Gamma_{0}^{(n)}|<
|\Gamma_{r}^{(n)}|< ...<|\Gamma_{r(m_1+m_2+1)}^{(n)}|$, then,
because $|\Gamma_{0}^{(n)}|\geq 2$, we obtain
$|\Gamma_{r(m_1+m_2+1)}^{(n)}|\geq {m_1+m_2+3}$. But on the other
hand $|\Gamma_{r(m_1+m_2+1)}^{(n)}|\leq
|B_{m_{1},m_{2}}^{(n)}|=m_1+m_2+3$. Consequently,
$|\Gamma_{r(m_1+m_2+1)}^{(n)}|=m_1+m_2+3$, and hence  the number of
elements of $\Gamma_{r(m_1+m_2+1)}^{(n)}$ coincides with the number
of elements of $B_{m_{1},m_{2}}^{(n)}$. In other words, every
element of the partition $\Gamma_{r(m_1+m_2+1)}^{(n)}$ contains only
one element of $B_{m_{1},m_{2}}^{(n)}$. Hence,
$|\Gamma_{s}^{(n)}|=m_1+m_2+3$ for all $s\geq r(m_1+m_2+1)$. We can
take $r(m_1+m_2+1)$ as $s_{0}$. It follows from the construction
that the number $s_{0}$ depends on $n$ but does not exceed
$r(m_1+m_2+1)$. The first and second claims of Theorem \ref{LG} are therefore
proved. Assume, that for every element $E_{s_{0}}^{(n)}$ of the
partition $\Gamma_{s_{0}}^{(n)}$ the relation
$$\prod\limits_{b_{1}:\overline{b}_{1}(n)\in
E_{s_{0}}^{(n)}}\sigma_{f_{1}}(b_{1})=
\prod\limits_{b_{2}:\overline{b}_{2}(n)\in
h(E_{s_{0}}^{(n)})}\sigma_{f_{2}}(b_{2}),$$
 holds. In this case,
$$
\prod\limits_{b_{1}\in B(f_1)}\sigma_{f_{1}}(b_{1})
=\prod\limits_{b_{2}\in B(f_2)}\sigma_{f_{2}}(b_{2}).
$$
in contradiction to the assumption in  Theorem \ref{LG}. This concludes the proof of  Theorem \ref{LG}.

\section{Jump coverings of the circle homeomorphisms with break points}

 We consider two $P$-homeomorphisms $f_{1}$ and $f_{2}$ with identical
irrational rotation number $\rho=\rho_{f_{1}}=\rho_{f_{2}}$.
Suppose that $f_{1}$ and $f_{2}$ has $m_{1}$ respectively $m_{2}$ break
points. Denote by $B(f_{i}),i=1,2$  the sets of all break points of
$f_{i}:$ $B(f_{1})=\{b_{1}^{(i)},\,\,1\leq i\leq m_1\}$ and
$B(f_{2})=\{b_{2}^{(i)},\,\,i=\overline{1,m_2}\}$

Next we introduce the notion of a "regular" cover of $B(f_{1})\cup
h^{-1}(B(f_{2}))$, that is the union of
 the set  of break points of $f_{1}$ and the $h$-preimage of the set  of break points
of $f_{2}$.

 Let $z_{i}\in S^{1},\,i=\overline{1,4},\,z_{1}\prec
z_{2}\prec z_{3}\prec z_{4}\prec z_{1}$ and let $r_{n}$ take values
in the set $\{q_{n-1},q_{n},q_{n-1}+q_{n}\}$. Suppose furthermore
that the interval $[z_{1},z_{4}]$ is $r_n$-small, i. e. the
intervals $\{f_{1}^{j}([z_{1},z_{4}]),\,0\leq j\leq r_{n}-1\}$, are
pairwise disjoint. Suppose that the system of intervals
$\{f_{1}^{j}([z_{1},z_{4}]),\,0\leq j\leq r_{n}-1\},$ covers the
elements of some non-empty subset $\widehat{B}_1\subset B(f_{1})$
with $\widehat{B}_1=\{b_{1}^{(i_s)},\, 1\leq s\leq p_1\}$. For every
element $b_{1}^{(i_s)}\in \widehat{B}_1$ there exists then a number
$l_{1}^{(i_{s})},\,0\leq l_{1}^{(i_{s})}\leq r_{n}-1$, such that
$\overline{b}_{1}^{(i_s)}(n)=f_{1}^{-l_{1}^{(i_{s})}}(b_{1}^{(i_{s})})\in
[z_{1},z_{4}]$. The point $\overline{b}_{1}^{(i_{s})}(n)$ is called
the $r_n$-preimage of the element $b_{1}^{(i_{s})}$ in $[z_1,z_4]$.
The set of $r_n$-preimages of elements of $\widehat{B}_1$
consists then of the elements
$\overline{b}_{1}^{(i_{1})}(n),\overline{b}_{1}^{(i_{2})}(n),...,\overline{b}_{1}^{(i_{p_{1}})}(n)$.

Define
$$
\xi_{f_{1}}(j):=\frac{\ell([f_{1}^{j}(z_{2}),f_{1}^{j}(z_{3})])}
{\ell([f_{1}^{j}(z_{1}),f_{1}^{j}(z_{2})])},\,\,
z_{f_{1}}^{(i_s)}(j):=
\frac{\ell([f_{1}^{j}(\overline{b}_{1}^{(i_s)}{(n)}),f_{1}^{j}(z_{2})])}
{\ell([f_{1}^{j}(z_{1}),f_{1}^{j}(z_{2})])},\,\, 1\leq s\leq
p_1,\,\,0\leq j<r_n.
$$
It follows easily from Lemma 2.2 that
\begin{equation}\label{evf}
e^{-v}\xi_{f_{1}}(0)\le \xi_{f_{1}}(j)\le e^v\xi_{f_{1}}(0),\;\,
e^{-v}z_{f_{1}}^{(i_{s})}(0)\le z_{f_{1}}^{(i_{s})}(j)\le
e^vz_{f_{1}}^{(i_{s})}(0),\, 1\leq s\leq p_1,
\end{equation}
and all $1\le j\le r_n-1$.
For the further discussion we introduce some definitions.
\begin{definition}\label{defCO}
Let $K>M\geq 1,\,\zeta\in (0,1),\,\delta>0$ be constant numbers, let
$n$ be a positive integer, and let $x_{0}\in S^{1}$. We say that a
triple of intervals
$([z_{1},z_{2}],[z_{2},z_{3}],[z_{3},z_{4}]),\,z_{i}\in
S^{1},\,i=1,2,3,4,$ covers
the break points in a subset $\widehat{B}_1$  "$(K,M,\delta,\zeta;x_{0})$-regularly" , if for some $r_{n}\in
\{q_{n-1},q_{n},q_{n}+q_{n-1}\}$ the following conditions hold:

1) $[z_{1},z_{4}]\subset (x_{0}-\delta,x_{0}+\delta),$ and the
system of intervals $\{f_{1}^{j}([z_{1},z_{4}]),\,0\leq j\leq
r_{n}-1\}$
covers every point in $\widehat{B}_1$ only once;\\

2) $\overline{b}_{1}^{(i_s)}(n)\in [z_1,z_2],\,\,1\leq s\leq p_1$;\\

3) $M\ell([z_{1},z_{2}])\leq \ell([z_{2},z_{3}])\leq
K\ell([z_{1},z_{2}])$,\,\,\,$K^{-1}\ell([z_{3},z_{4}])\leq
\ell([z_{2},z_{3}])\leq K\ell([z_{3},z_{4}])$;\\

4) The lengths of the intervals
$f_{1}^{r_{n}}([z_{1},z_{2}]),\,f_{1}^{r_{n}}([z_{2},z_{3}])$ and
$f_{1}^{r_{n}}([z_{3},z_{4}])$ are pairwise $K$-comparable;\\

5)
$\max\{\ell([f_{1}^{r_{n}}(z_{i}),x_{0}]),\ell([z_{i},x_{0}]),i=\overline{1,4}\}
\leq K\ell([z_{1},z_{2}])$;\\

6) $\max\limits_{1\leq s\leq p_{1}}\{z_{f_{1}}^{(i_{s})}(0)\}<\zeta$.\\
\end{definition}

\begin{definition}\label{defNC} Let $\widehat{B}_1$ and $\widehat{B}_2$
be subsets of the break points of the homeomorphisms
  $f_{1}$ and $f_{2}$, respectively. The subsets $\widehat{B}_1$ and
$\widehat{B}_2$ are said to be \textbf{"not jump-coinciding"}, or for short \textbf{"not coinciding"},if
$$
\prod\limits_{b_{1}\in \widehat{B}_1}\sigma_{f_{1}}(b_{1})\not
=\prod\limits_{b_{2}\in \widehat{B}_2}\sigma_{f_{2}}(b_{2}).
$$
Otherwise we call them \textbf{"coinciding"}
\end{definition}

It is clear, that if $\widehat{B}_1$ and $\widehat{B}_2$ are "not
jump-coinciding" subsets, then one of subsets  $\widehat{B}_{1}$ and
$\widehat{B}_{2}$ is non-empty.  For instance, if  $\widehat{B}_1$
is empty, then we put $\prod\limits_{b_{1 }\in
\widehat{B}_1}\sigma_{f_{1}}(b_{1}):=1$.

\begin{definition}\label{defCO1} Let $\widehat{B}_1$ and $\widehat{B}_2$ be
subsets of the break points of $f_{1}$ and $f_{2}$ respectively. We say
that the triples of intervals $([z_s,z_{s+1}],\,\,s=1,2,3)$ and
$([h(z_s),h(z_{s+1})],\,\,s=1,2,3)$
cover the subsets
$\widehat{B}_1$ and $\widehat{B}_2$ "$(K,M,\delta,\zeta;x_{0},h(x_0))$- regularly" with $r_{n}\in
\{q_{n-1},q_{n},q_{n}+q_{n-1}\}$, if

$\bullet$\,\, the triples of intervals $([z_s,z_{s+1}],\,\,s=1,2,3)$
and $([h(z_s),h(z_{s+1})],\,\,s=1,2,3)$
 cover the points of
$\widehat{B}_1$ respectively $\widehat{B}_2$
 "$(K,M,\delta,\zeta;x_{0})$- regularly" respectively
"$(K,M,\delta,\zeta;h(x_0))$-regularly" if $\widehat{B}_1\neq {\O},\,\,\widehat{B}_2\neq {\O}$;\\

$\bullet$\,\, the triple of intervals $([z_s,z_{s+1}],\,\,s=1,2,3)$
 covers the points of
$\widehat{B}_1$ "$(K,M,\delta,\zeta;x_{0})$-regularly" if $\widehat{B}_1\neq {\O},\,\,\widehat{B}_2=
{\O}$;\\

$\bullet$\,\, the triple intervals
$([h(z_s),h(z_{s+1})],\,\,s=1,2,3)$
 covers the points of $\widehat{B}_2$ "$(K,M,\delta,\zeta;\\h(x_0))$-regularly" if $\widehat{B}_1={\O},\,\,\widehat{B}_2\neq
{\O}$.
\end{definition}

Next we formulate our main Theorem on the covering of intervals which plays
a key role in the proof of Theorem \ref{DMS}.

\begin{theorem}\label{CO} Suppose that the homeomorphisms $f_{1}$ and $f_{2}$
satisfies the conditions of Theorem \ref{DMS}. Let
$Dh(x_{0})=\omega_0>0$ for some $x_0\in S^{1}$ and
let $M\geq 1,\,\zeta,\,\delta\in (0,1)$ be constants. Then there
exist a constant $K=K(f_{1},f_{2},M,\zeta)>M$ and for any
sufficiently large $n$  "not jump-coinciding" subsets
$\widehat{B}_1$ and $\widehat{B}_2$, points $z_{i}\in S^{1},\,1\leq
i\leq 4$ with $z_{1}\prec z_{2}\prec z_{3}\prec z_{4}\prec z_{1}$
and a number $r_{n}=r_{n}(z_{1},z_{2},z_{3},z_{4})\in
\{q_{n-1},q_{n},q_{n}+q_{n-1}\},$ such that  the triples of
intervals $([z_s,z_{s+1}],\,\,s=1,2,3)$ and
$([h(z_s),h(z_{s+1})],\,\,s=1,2,3)$
 cover the  subsets
$\widehat{B}_1$ and $\widehat{B}_2$ "$(K,M,\delta,\zeta;x_{0},h(x_0))$-regularly" for $r_{n}$.
\end{theorem}

\textbf{Proof of Theorem \ref{CO}.}   Suppose that the
homeomorphisms $f_{1}$ and $f_{2}$ satisfy the conditions of Theorem
\ref{DMS}. Let $B(f_{1})=\{b_{1}^{(i)},\,\,1\leq i\leq m_1\}$ and
$B(f_{2})=\{b_{2}^{(j)},\,\,1\leq j\leq m_2\}$ be the sets of break
points of $f_{1}$ and $f_{2}$ respectively. By assumption
$Dh(x_{0})=\omega_0>0$ for some $x_{0}\in S^{1}.$
Consider the dynamical partition $\eta_{n}(x_0)$ of
the point $x_{0}$ under $f_{1}.$ Suppose $n$ to be  odd. Let
$B_{m_{1},m_{2}}^{(n)}$ and $d_{n}$ be defined as in (\ref{BMM}).

Define the number $m_0$ by using the constants
$M>1,\,\,\zeta\in(0,1)$ and the total variation $v_i$ of the
functions $\ln Df_i,\,\,i=1,2$ as follows:
\begin{equation}{\label{m_0}}
m_{0}:=\max\{m_1+m_2+4,\,[{M}{\zeta}^{-1}]+1,\,[e^{v_1}]+1,\,[e^{v_2}]+1\},
\end{equation}
where $[\cdot]$ denotes the integer part, and consider the partition
$D_{l}^{(n)}$ of the interval $[x_{q_{n}},x_{q_{n-1}}].$ It is
sufficient to use the assertion of Theorem \ref{LG} with $r=9$ and
set $s_{0}=s_{0}(9,n)$. By this assertion  there exists at least one
element $\widetilde{E}_{s_{0}}^{(n)}\in \Gamma_{s_0}^{(n)}$ such
that
$$\prod\limits_{b_{1}:\overline{b}_{1}(n)\in
\widetilde{E}_{s_{0}}^{(n)}}\sigma_{f_{1}}(b_{1})\neq
\prod\limits_{b_{2}:\overline{b}_{2}(n)\in
h(\widetilde{E}_{s_{0}}^{(n)})}\sigma_{f_{2}}(b_{2}).$$

Set  $\widehat{B}_1:=\{b_{1}^{(i)}:\overline{b}_{1}^{(i)}(n)\in
\widetilde{E}_{s_{0}}^{(n)}\}$ and
$\widehat{B}_2:=\{b_{2}^{(j)}:\overline{b}_{2}^{(j)}(n)\in
h(\widetilde{E}_{s_{0}}^{(n)})\}.$ Then the
following cases are possible:
 $\widehat{B}_1\neq {\O},\widehat{B}_2\neq {\O}$ or
$\widehat{B}_1\neq {\O},\widehat{B}_2={\O}$ or $\widehat{B}_1=
{\O},\widehat{B}_2\neq {\O}$. If $\widehat{B}_1\neq
{\O}$ and $\widehat{B}_2\neq {\O}$ we'll construct a triple of
regular covering intervals $[z_{s},z_{s+1}],\,\,s=1,2,3$, in the other
cases the construction of regular covering intervals is analogous.

Let $\widehat{B}_1\neq {\O}$ and $\widehat{B}_2\neq {\O}$. Then three
cases are possible for the set $\widetilde{E}_{s_{0}}^{(n)}$.

$(c_1)$ $\widetilde{E}_{s_{0}}^{(n)}$ does not contain any elements
of $\{x_{q_{n}},\,x_{0},\,x_{q_{n-1}}\}$;

$(c_2)$  $\widetilde{E}_{s_{0}}^{(n)}$ contains only one element of
the set $\{x_{q_{n}},\,x_{0},\,x_{q_{n-1}}\}$;

$(c_3)$ $\widetilde{E}_{s_{0}}^{(n)}$ contains the elements
$x_{q_{n}},\,x_{0}$ of the set
$\{x_{q_{n}},\,x_{0},\,x_{q_{n-1}}\}$.\newline  The case
 $\{x_{q_{n-1}},\,x_{0}\}\in \widetilde{E}_{s_{0}}^{(n)}$ turns out to be impossible.\newline
 We prove the assertion of
the theorem in each of the cases separately.

$(c_1)$.  Let
$\widetilde{E}_{s_{0}}^{(n)}\cap{\{x_{q_{n}},\,x_{0},\,x_{q_{n-1}}\}}=\O$.
Then either $\widetilde{E}_{s_{0}}^{(n)}\subset{(x_{q_{n}},x_{0})},$
or $\widetilde{E}_{s_{0}}^{(n)}\subset{(x_{0},x_{q_{n-1}})}$.
Suppose for definiteness that
$\widetilde{E}_{s_{0}}^{(n)}\subset{(x_{0},x_{q_{n-1}})}$. The case
$\widetilde{E}_{s_{0}}^{(n)}\subset{(x_{q_{n}},x_{0})}$ can be treated
in a similar way.

One can deduce from the assertion of Theorem \ref{LG} that the
subset $\widetilde{E}_{s_{0}}^{(n)}$ is covered by one or two
intervals of the partition $D_{s_0}^{(n)}$. The
union of the intervals of the partition $D_{s_0}^{(n)}$ which cover
$\widetilde{E}_{s_{0}}^{(n)}$ is denoted by $I_{s_0}^{(n)}$. In the
same way we can define intervals $I_{s_0+p}^{(n)}$ for $0<p<9$.
Clearly $I_{s_{0}}^{(n)}\supset I_{s_{0}+1}^{(n)}...\supset
I_{s_{0}+9}^{(n)}$. It follows from the assertion of Theorem
\ref{LG} that the interval $I_{s_0}^{(n)}$ is adjacent on the left
and right to two $s_{0}^{(n)}$-empty intervals of $D_{s_0}^{(n)}$
contained in the interval $(x_{0},x_{q_{n-1}})$. These two intervals
are denoted by $L_{s_{0}}^{(n)}$ and $R_{s_{0}}^{(n)}$,
respectively.

We now define the points $z_{i},\,1\leq i\leq 4$, as follows:
$$
z_{2}=\max \{y: y\in\widetilde{E}_{s_{0}}^{(n)}\},\,\,\,z_{1}=z_{2}-d_{n}
m_{0}^{-(s_{0}+7)},\,\,\,z_{3}= z_{2}+d_{n} m_{0}^{-(s_{0}+6)},
$$
\begin{equation}\label{Z1} z_{4}=z_{2}+d_{n} m_{0}^{-(s_{0}+6)}+d_{n} m_{0}^{-(s_{0}+7)}.
\end{equation}
We now verify that the triples of intervals
$[z_{s},z_{s+1}],\,\,s=1,2,3$ and
$[h(z_{s}),h(z_{s+1})],\,\,s=1,2,3$ satisfy the  conditions of
Definition \ref{defCO}. The length of the interval
$I_{s_0+9}^{(n)}$ covering the subset
$\widetilde{E}_{s_{0}}^{(n)}$ does not exceed $2d_{n}
m_{0}^{-(s_{0}+9)}$, and the lengths of the intervals
$L_{s_{0}}^{(n)}$ and $R_{s_{0}}^{(n)}$ adjacent to $I_{s_0}^{(n)}$
are equal to $d_{n}\cdot m_{0}^{-s_{0}}$. Using the definition of
the points $z_{i},0\leq i\leq 4$  we obtain
$\ell([z_{1},z_{2}])=m_{0}^{-7}\ell(L_{s_{0}}^{(n)}),\,\ell([z_{2},z_{4}])=
(m_{0}+1)m_{0}^{-7}\ell(R_{s_{0}}^{(n)})$. Hence,
$[z_{1},z_{4}]\subset{L_{s_{0}}^{(n)}\cup I_{s_{0}}^{(n)}\cup
R_{s_{0}}^{(n)}}\subset{(x_{0},x_{q_{n-1}})}$. Since the interval
$[x_{0},x_{q_{n-1}}]$ is $q_n$-small, the intervals
$\{f_{1}^{j}([z_{1},z_{4}]),\,0\leq j\leq q_{n}-1 \}$ are pairwise
disjoint and cover each point of  $\widehat{B}_1$ only once. One can
easily verify that the intervals $[z_{s},z_{s+1}],\,\,s=1,2,3$
satisfy condition 2) of Definition \ref{defCO}. By Denjoy's
inequality the intervals $[x_{q_n},x_{q_{n-1}}]$ and $[x_0,
x_{q_{n-1}}]$ are $1+e^{v_1}$-comparable. Hence, using the fact that
$\ell([z_{1},z_{4}])=(m_{0}+2)m_{0}^{-(s_{0}+7)}d_{n}$, we obtain
that $[z_{1},z_{4}]$ and $[x_{0},x_{q_{n-1}}]$ are
$(m_{0}+2)^{-1}m_{0}^{(s_{0}+7)}(1+e^{v_{1}})^{-1}$-
comparable. Set
$K:=\max\{m_{0}e^{2(v_{1}+v_{2})},m_{0}^{s_{0}+7},\,n=1,2,...\},$
where $s_{0}=s_{0}(9,n)$. By the assertion of Theorem \ref{LG} we
have $s_{0}=s_{0}(9,n)\leq 9(m_1+m_2+1),\,\forall{n}\in N$.
Consequently, $K=\max\{m_{0}e^{2(v_1+v_2)},m_{0}^{9(m_1+m_2)+16}\}$.
Taking into account that $m_{0}>M$ we conclude that the intervals
$[z_{s},z_{s+1}],\,\,s=1,2,3$, satisfy condition 3) of Definition
\ref{defCO} with the constant $K$. By Denjoy's inequality the
intervals $[z_{s},z_{s+1}]$ and
$[f_{1}^{q_{n}}(z_{s}),f_{1}^{q_{n}}(z_{s+1})]$ are
$e^{v_1}$-comparable for every $s=1,2,3$. Since
$\ell([z_{s},z_{s+1}])=d_{n}\cdot m_{0}^{-(s_{0}+7)},\,s=1,3$ and
$\ell([z_{2},z_{3}])=d_{n}\cdot m_{0}^{-(s_{0}+6)}$, it follows that
the intervals
$[f_{1}^{q_{n}}(z_{s}),f_{1}^{q_{n}}(z_{s+1})],\,s=1,2,3$, satisfy
condition 4) of Definition \ref{defCO} with the constant $K$.

Obviously,
$$
\max\limits_{1\leq i\leq 4}\ell([z_{i},x_{0}]),\,\max\limits_{1\leq
i\leq 4}\ell([f_{1}^{q_{n}}(z_{i}),x_{0}])\leq
\ell([x_{q_{n}},x_{q_{n-1}}])=d_{n}.
$$
It follows from the explicit form of the length of the interval
$\ell([z_{1},z_{2}])$ that
$d_{n}=m_{0}^{s_{0}+7}\\\ell([z_{1},z_{2}])$. Hence it is evident
that $[z_{s},z_{s+1}],\,s=1,2,3,$ satisfy condition 5) of Definition
\ref{defCO} with the constant $K$. We now verify that the triple of
intervals $[z_{s},z_{s+1}],\,s=1,2,3,$ satisfy condition 6) of
Definition \ref{defCO}. It follows from the definition of the points
$z_{s},\,s=1,2,3,4$, that
$$
\max\limits_{i:\overline{b}_{1}^{(i)}(n)\in
\widetilde{E}_{s_{0}}^{(n)}}
\frac{\ell([\overline{b}_{1}^{(i)}(n),z_{2}])}{\ell([z_{1},z_{2}])}\leq
\frac{d_{n}\cdot m_{0}^{-(s_{0}+8)}}{d_{n}\cdot
m_{0}^{-(s_{0}+7)}}=m_{0}^{-1}<\zeta.
$$

Next we show that the intervals $[h(z_{s}),h(z_{s+1})],\,\,s=1,2,3,$
 define a regular cover of the elements of  the set $\widehat{B}_2$.
The definition of the points $z_{i},\,1\leq i \leq 4$ implies that
$[h(z_{1}),h(z_{4})]\subset[h(x_0),f_{2}^{q_{n-1}}(h(x_0))]\subset(h(x_0)-\delta,h(x_0)+\delta)$
for sufficiently large $n$ and the system of intervals
$\{f_{2}^{j}([h(z_{1}),h(z_{4})])$,\,$0\leq j<q_n\}$ are pairwise
disjoint and cover each  point of the subset $\widehat{B}_2$ only
once. Since $\widetilde{E}_{s_{0}}^{(n)}\subset [z_{1},z_{4}]$, the
elements of the set $\widehat{B}_2$ are covered by
$f_{2}$- iterations of the interval
$[h(z_{1}),h(z_{2})]$. By the assumption of  Theorem \ref{CO}
$Dh(x_0)=\omega_0>0$ at the point $x_0\in S^{1}$.

 Let $H(x)$ be the  lift of $h.$ By the definition
of the derivative, for any $\varepsilon>0$ there exists
$\delta_{1}=\delta_{1}(x_0,\varepsilon)>0$ such that
for all $x\in (x_0-\delta_{1},x_0+\delta_{1})$ the inequality
\begin{equation}\label{OMEGA}
\omega_0-\varepsilon<\frac{H(x)-H(x_0)}{x-x_0}<\omega_0+\varepsilon.
\end{equation}
holds. For $x=\hat{z}_i,\,\,1\leq i\leq 4$, it follows from (\ref{OMEGA}) that
\begin{equation}\label{OMEGA1}
(\omega_0-\varepsilon)(x_{0}-\hat{z}_i)<H(x_0)-H(\hat{z}_i)<(\omega_0+\varepsilon)(x_{0}-\hat{z}_i),
\end{equation}
from which  one can easily derive the inequalities
$$
\omega_0-\varepsilon\frac{(x_0-\hat{z}_{i+1})+(x_0
-\hat{z}_{i})}{\hat{z}_{i+1}-\hat{z}_{i}}<
\frac{H(\hat{z}_{i+1})-H(\hat{z}_{i})}{\hat{z}_{i+1}-\hat{z}_{i}}<
$$
\begin{equation}\label{OMEGA2}
<\omega_0+\varepsilon\frac{(x_0-\hat{z}_{i+1})+(x_0-\hat{z}_{i})}
{\hat{z}_{i+1}-\hat{z}_{i}},\,\,i=1,2,3,
\end{equation}
where $(\hat{z}_{1},\hat{z}_{2},\hat{z}_{3},\hat{z}_{4})$ is the
lifted vector of $(z_{1},z_{2},z_{3},z_{4})$. Using the definition
of the points $z_{i},\,1\leq i\leq 4,$  we obtain
\begin{equation}\label{OMEGA3}
\max\limits_{1\leq i\leq
3\{\frac{x_0-\hat{z}_{i+1}}{\hat{z}_{i+1}-\hat{z}_{i}}
,\,\frac{x_0-\hat{z}_{i}}{\hat{z}_{i+1}-\hat{z}_{i}}\}<K.}
\end{equation}
Using the definition of $m_0$ and the bounds
(\ref{OMEGA2}),(\ref{OMEGA3}) it can be easily shown that the
intervals $[h(z_{s}),h(z_{s+1})],\,\,s=1,2,3$ satisfy the other
conditions of Definition \ref{defCO} with the same constant $K$.

$(c_2)$. We denote by
$E_{s_{0}}^{(n)}(x_{q_{n}}),\,E_{s_{0}}^{(n)}(x_{0})$ and
$E_{s_{0}}^{(n)}(x_{q_{n-1}})$ the elements of the partition
$\Gamma_{s_{0}}^{(n)}$ containing just one of the points
$x_{q_{n}}$, $x_{0}$ or $x_{q_{n-1}}$, and which are pairwise
disjoint. By assumption, $\widetilde{E}_{s_{0}}^{(n)}$ coincides
with  one of them and these
subsets are pairwise disjoint. Suppose, that for any
element $E_{s_{0}}^{(n)}\in \Gamma_{s_0}^{(n)}$,  which does not contain the points $x_{q_n}$,
$x_0$, $x_{q_{n-1}}$ the following relation hold:
$$
\prod\limits_{b_{1}:\overline{b}_{1}(n)\in
E_{s_{0}}^{(n)}}\sigma_{f_{1}}(b_{1})=
\prod\limits_{b_{2}:\overline{b}_{2}(n)\in
h(E_{s_{0}}^{(n)})}\sigma_{f_{2}}(b_{2}).
$$
Otherwise we arrive again at case $(c_{1})$. Then the sets
$\widetilde{B}_1=\{b_{1}^{(i)}:\overline{b}_{1}^{(i)}(n)\in
E_{s_{0}}^{(n)}(x_{q_{n}})\cup E_{s_{0}}^{(n)}(x_{0})\cup
E_{s_{0}}^{(n)}(x_{q_{n-1}})\}$ and
$\widetilde{B}_2=\{b_{2}^{(i)}:h^{-1}(\overline{b}_{2}^{(i)}(n))\in
E_{s_{0}}^{(n)}(x_{q_{n}})\cup E_{s_{0}}^{(n)}(x_{0})\cup
E_{s_{0}}^{(n)}(x_{q_{n-1}})\}$ are "not coinciding".

Denote by $I_{s_0+p}^{(n)}(x_{q_n})$,
$I_{s_0+p}^{(n)}(x_0)$ and $I_{s_0+p}^{(n)}(x_{q_{n-1}})$ be the
intervals of  partition $D_{s_{0}+p}^{(n)},$ $1\leq p\leq 9,$
covering $E_{s_{0}}^{(n)}(x_{q_{n}})$, $E_{s_{0}}^{(n)}(x_{0})$ and
$E_{s_{0}}^{(n)}(x_{q_{n-1}})$, respectively. By Theorem \ref{LG}
the interval $I_{s_0}^{(n)}(x_0)$ is adjacent on both sides to
$s_{0}^{(n)}$-empty intervals of length $d_{n}m_{0}^{-s_0}$.
Consider then the subset
$E_{s_{0}}^{(n)}(x_{q_{n}},x_{0},x_{q_{n-1}}):=f_{1}^{-q_{n}}(E_{s_{0}}^{(n)}(x_{q_{n}}))\cup
E_{s_{0}}^{(n)}(x_{0})\cup
f_{1}^{-q_{n-1}}(E_{s_{0}}^{(n)}(x_{q_{n-1}}))$. From Denjoy's
inequality we get
$\ell(f_{1}^{-q_n}(I_{s_{0}+9}^{(n)}(x_{q_{n}})))$,\,
$\ell(f_{1}^{-q_{n-1}}(I_{s_{0}+9}^{(n)}(x_{q_{n-1}})))\leq
e^{v_1}d_{n}m_{0}^{-(s_{0}+9)}$. Hence,
$E_{s_{0}}^{(n)}(x_{q_{n}},x_{0},x_{q_{n-1}})
\subset(x_{0}-2e^{v_{1}}d_{n}m_{0}^{-(s_{0}+9)},x_{0}+
2e^{v_{1}}d_{n}m_{0}^{-(s_{0}+9)})$.

We can define now the points $z_i,\, 1 \leq i\leq 4$, as follows
$$
z_{2}=\max
E_{s_{0}}^{(n)}(x_{q_{n}},x_{0},x_{q_{n-1}}),\,z_{1}=z_{2}-d_{n}m_{0}^{-(s_{0}+6)},\,
z_{3}=z_{2}+d_{n}m_{0}^{-(s_{0}+3)},
$$
\begin{equation}
\label{Z2} z_{4}=z_{2}+d_{n}m_{0}^{-(s_{0}+3)}+d_{n}m_{0}^{-(s_{0}+6)}.
\end{equation}

We claim, the triple of intervals $[z_{s},z_{s+1}],\,s=1,2,3$
covers the point of
the subset $\widetilde{B}_1$  "$(K,M,\delta,\zeta;x_{0})$-regularly" with $r_{n}=q_{n}+q_{n-1}$.

We  verify only that the system of intervals
$\{f_{1}^{i}([z_{1},z_{4}]),\,0\leq i\leq q_{n}+q_{n-1}-1\}$ covers
each point of the subset
$E_{s_{0}}^{(n)}(x_{q_{n}},x_{0},x_{q_{n-1}})$ only once. The other
conditions in Definition \ref{defCO} concerning the lengths of the
intervals can be verified as in the case
$(c_{1})$ by simple calculations . We divide the interval $[z_{1},z_{4}]$ up into
$[z_{1},z_{4}]=[z_{1},x_{0}]\cup (x_{0},z_{4}]$. We claim, the
intervals $f_{1}^{i}([z_{1},x_{0}]),\,0\leq i\leq q_{n}+q_{n-1}-1$,
cover the break points of $f_{1}$ with $q_{n-1}$-preimages in
$E_{s_{0}}^{(n)}(x_{0})\cap [z_{1},x_{0}]$ and with $q_{n}$-preimages in
$E_{s_{0}}^{(n)}(x_{q_{n-1}})$ only once. Since
$[z_{1},x_{0}]\subset[x_{q_{n}},x_{0}]$, the intervals
$f_{1}^{i}([z_{1},x_{0}]),\,0\leq i\leq q_{n-1}-1$, cover each break
point of $f_{1}$ with $q_{n}$-preimage in $[z_{1},x_{0}]\cap
E_{s_{0}}^{(n)}(x_{0})$ only once. It follows from the assertion of
Theorem \ref{LG} that there is an $s_{0}$-empty interval to the left
of $I_{s_{0}}^{(n)}(x_{q_{n}})$. By Denjoy's inequality the length
of the interval
$f_{1}^{q_{n-1}}([z_{1},x_{0}])=[f_{1}^{q_{n-1}}(z_{1}),x_{q_{n-1}}]$
is at most $e^{v_1}\ell([z_{1},x_{0}])$. It is easy to verify that
this number is less than the sum of the lengths of the intervals
$I_{s_{0}}^{(n)}(x_{q_{n-1}})$ and the adjacent $s_{0}^{(n)}$-empty
one. In other words, $[f_{1}^{q_{n-1}}(z_{1}),x_{q_{n-1}}]$ is
covered by the interval $I_{s_{0}}^{(n)}(x_{q_{n-1}})$ and an
$s_{0}^{(n)}$-empty adjacent interval. Clearly, the subset of the
$q_{n}$- preimages of the break points of $f_{1}$ contained in
$[f_{1}^{q_{n-1}}(z_{1}),x_{q_{n-1}}]$ coincides with
$E_{s_{0}}^{(n)}(x_{q_{n-1}})$. Since
$[f_{1}^{q_{n-1}}(z_{1}),x_{q_{n-1}}]\subset[x_{0},x_{q_{n-1}}]$,
the intervals
$f_{1}^{i}([f_{1}^{q_{n-1}}(z_{1}),x_{q_{n-1}}]),\,0\leq i\leq
q_{n}-1$, cover the break points with $q_{n}$-preimages in
$[f_{1}^{q_{n-1}}(z_{1}),x_{q_{n-1}}]$ only once.

It can be shown in a similar way that the intervals
$f_{1}^{i}((x_{0},z_{4}]),\,0\leq i\leq q_{n}+q_{n-1}-1$, cover the
break points of $f_{1}$ with $q_{n}$-preimages in
$E_{s_{0}}^{(n)}(x_{0})\cap (x_{0},z_{4}]$ respectively $q_{n-1}$-preimages
in $E_{s_{0}}^{(n)}(x_{q_{n}})$ only once. Remember, $\varepsilon$
to be an arbitrary positive number. Using this and the bounds
(\ref{OMEGA2}),(\ref{OMEGA3}) it can be proved that the intervals
$[h(z_{s}),h(z_{s+1})],\,\,s=1,2,3$ satisfy all the conditions of
Definition \ref{defCO} with the same constant $K$.

 $(c_3)$. First we
show that the subset $\{x_0,x_{q_{n-1}}\}$ cannot be part of
$\widetilde{E}_{s_{0}}^{(n)}$. Suppose on the contrary $\{x_0,x_{q_{n-1}}\}\subset\widetilde{E}_{s_{0}}^{(n)}$ . Then
$\widetilde{E}_{s_{0}}^{(n)}$ is covered by the interval
$I_{s_{0}+9}^{(n)}$ of the partition $D_{s_{0}+9}^{(n)}$.
Consequently $\Delta_{0}^{(n-1)}\subset
I_{s_{0}+9}^{(n)}$. Clearly, $\ell(\Delta_{0}^{(n-1)})\leq
m_{0}^{-(s_{0}+9)}d_{n}$. Hence, $\ell(\Delta_{0}^{(n)})\geq
(m_{0}^{(s_{0}+9)}-1)\ell(\Delta_{0}^{(n-1)})$. From
Denjoy's inequality we get  $\ell(\Delta_{0}^{(n)})\leq
e^{v}\ell(\Delta_{0}^{(n-1)})$. Using the definition of $m_0$, one
can easily show that $m_{0}^{(s_{0}+9)}-1>e^{v_{1}}$. Consequently,
$\ell(\Delta_{0}^{(n)})>\ell(\Delta_{0}^{(n)})$,
a contradiction.

By assumption $(c_3)$ $x_{q_{n}},x_{0}\in
\widetilde{E}_{s_{0}}^{(n)}$. The definition of the subsets
$E_{s_{0}}^{(n)}(x_{q_{n}})$ and $E_{s_{0}}^{(n)}(x_{0})$ imply that
$\widetilde{E}_{s_{0}}^{(n)}=E_{s_{0}}^{(n)}(x_{0})=E_{s_{0}}^{(n)}(x_{q_{n}})$.

We divide the set
$E_{s_{0}}^{(n)}(x_{q_{n-1}})\setminus\{x_{q_{n-1}}\}$  up into two subsets
$A_{1}$ and $A_{2}$:

1)$A_{1}=\{\overline{b}_{1},h^{-1}(\overline{b}_{2}):\overline{b}_{1},h^{-1}(\overline{b}_{2})\in
E_{s_{0}}^{(n)}(x_{q_{n-1}}),\,\,f_{1}^{j}(\overline{b}_{1})$ and
$f_{2}^{s}(\overline{b}_{2})$ are break points of $f_1$ and $f_2$,
respectively for some $j,s$ with $ q_n - q_{n-1} \leq j,s <
q_{n-1}$\};

2)$A_{2}=\{\overline{b}_{1},h^{-1}(\overline{b}_{2}):\overline{b}_{1},h^{-1}(\overline{b}_{2})\in
E_{s_{0}}^{(n)}(x_{q_{n-1}}),$ $f_{1}^{i}(\overline{b}_{1})$ and
$f_{2}^{k}(\overline{b}_{2})$ are break points of $f_{1}$ and
$f_{2}$, respectively for some $i,k,$ with $0< i,k \leq q_n -
q_{n-1}-1\}$.

Let $I_{s_{0}+p}^{(n)}(x_{0})$ and
$I_{s_{0}+p}^{(n)}(x_{q_{n-1}}),\,0\leq p\leq 9,$  be the intervals
of partition $D_{s_0+p}^{(n)}$ covering $E_{s_{0}}^{(n)}(x_{0})$ and
$E_{s_{0}}^{(n)}(x_{q_{n-1}})$, respectively. By the definition of
the partition $D_{l}^{(n)}$ we have
\begin{equation}\label{ELLI}
\ell(I_{s_{0}+p}^{(n)}(x_{0}))=\ell(I_{s_{0}+p}^{(n)}(x_{q_{n-1}}))=d_{n}m_{0}^{-(s_{0}+p)},\,0\leq
p\leq 9.
\end{equation}
From Denjoy's lemma it folows that
\begin{equation}\label{ELLI1}
\ell(f^{q_{n}-q_{n-1}}(I_{s_{0}+9}^{(n)}(x_{q_{n-1}}))),\,\ell(f^{-q_{n-1}}
(I_{s_{0}+9}^{(n)}(x_{q_{n-1}})))<d_{n}m_{0}^{-(s_{0}+7)}.
\end{equation}

We set
$E_{s_{0}}^{(n)}(x_{0},x_{q_{n-1}}):=E_{s_{0}}^{(n)}(x_{0})\cup
f_{1}^{q_{n}-q_{n-1}}(A_{1})\cup f_{1}^{-q_{n-1}}(A_{2})$.
Obviously, $E_{s_{0}}^{(n)}(x_{0},x_{q_{n-1}})\\\subset
H_{s_{0}+9}^{(n)}:= I_{s_{0}+9}^{(n)}(x_{0})\cup
f_{1}^{q_{n}-q_{n-1}}(I_{s_{0}+9}^{(n)}(x_{q_{n-1}}))\cup
f_{1}^{-q_{n-1}}(I_{s_{0}+9}^{(n)}(x_{q_{n-1}})) $. This fact
together with relations (\ref{ELLI}) and (\ref{ELLI1}), implies that
$\ell(H_{s_{0}+9}^{(n)})<2d_{n}m_{0}^{-(s_{0}+7)}$. We can now
define points $z_{i},\,i=1,2,3,4$ as follows:
$$z_{2}=\max
E_{s_{0}}^{(n)}(x_{0},x_{q_{n-1}}),\,z_{1}=z_{2}-2d_{n}m_{0}^{-(s_{0}+6)},\,
z_{3}=z_{2}+2d_{n}m_{0}^{-(s_{0}+4)},
$$
\begin{equation}\label{Z3}
z_{4}=z_{2}+d_{n}m_{0}^{-(s_{0}+4)}+d_{n}m_{0}^{-(s_{0}+6)}.
\end{equation}

We claim that the triples of intervals $([z_{s},z_{s+1}],s=1,2,3)$
and $([h(z_{2}),h(z_{3})],s=1,2,3)$
 cover the points of
the "not coinciding" sets
$\widetilde{B}_1=\{b_{1}^{(i)}:\overline{b}_{1}^{(i)}(n)\in
E_{s_{0}}^{(n)}(x_{0})\cup E_{s_{0}}^{(n)}(x_{q_{n-1}})\}$ and
$\widetilde{B}_2=\{b_{2}^{(i)}:h^{-1}(\overline{b}_{2}^{(i)}(n))\in
E_{s_{0}}^{(n)}(x_{0})\cup E_{s_{0}}^{(n)}(x_{q_{n-1}})\}$ "$(K,M,\delta,\zeta;x_{0},h(x_0))$- regularly" with
$r_n=q_n$.

We shall only show that condition 1) in Definition \ref{defCO}
holds, the other conditions can easily be verified as in case
$(c_1)$. We again divide the interval $[z_{1},z_{4}]$ up into
$[z_{1},z_{4}]=[z_{1},x_{0}]\cup [x_{0},z_{4}]$ and decompose the
system of intervals $\{f_{1}^{i}([z_{1},x_{0}]),\,0\leq i\leq
q_{n}-1\}$ into two subsystems: $\{f_{1}^{i}([z_{1},x_{0}]),\,0\leq
i\leq q_{n-1}-1\}$ and $\{f_{1}^{i}([z_{1},x_{0}]),\,q_{n-1}\leq
i\leq q_{n}-1\}$. Clearly, the first subsystem covers  those break
points of $f_1$ whose $q_{n-1}$- preimages are contained in the
interval $[x_{q_n}, x_0]$, as well as those whose $q_n$-pre- images
form the subset $A_{1}$, only once. Note, that the break points of
$f_1$ with $q_n$-preimages in $A_2$ are covered by the second system
of intervals only once. Consider then the system of intervals
$\{f_{1}^{i}([x_{0},z_{4}]),\,0\leq i\leq q_{n}-1\}$. It follows
from the definition of $z_{4}$ that $[x_{0},z_{4}]\subset
[x_{0},x_{q_{n-1}}]$. Then two cases are possible for the point
$z_{2}$: either $z_{2}>x_{0}$ or $z_{2}\leq x_{0}$. If
$z_{2}>x_{0}$, then the break points with $q_n$-preimages contained
in $[x_0,z_2]$ are covered by the system of intervals
$\{f_{1}^{i}([x_{0},z_{4}]),\,0\leq i\leq q_{n}-1\}$ only once. In
the case $z_2<x_0$, the interval $[x_0,z_4]$ does not contain any
$q_n$-preimages of break points of $f_1$. Notice that the orbits
$\{f_{1}^{k}(x_0),\,\,k\in \mathbb{Z}\}$ and
$\{f_{2}^{k}(h(x_0)),\,\,k\in \mathbb{Z}\},$ of any point
 $ x_0\in S^1$ has the same order on the circle. This together with (\ref{OMEGA2})
and (\ref{OMEGA3}) implies that the intervals
$[h(z_{s}),h(z_{s+1})],\,\,s=1,2,3$ also satisfy the conditions of
Definition \ref{defCO} with the same constant $K$. Hence the theorem is
completely proved.
\section{Proof of the Theorem \ref{DMS}}

Consider two copies of the circle $S^1$, and homeomorphisms  $f_{1}$
and $f_{2}$ with $m_1,m_2\geq 2$
break points respectively and identical irrational rotation number
$\rho$. Denote by
 $B(f_{1})=\{b_{1}^{(i)},1\leq i \leq m_{1}\}$ and
$B(f_{2})=\{b_{2}^{(j)},1\leq j\leq m_{2}\}$ the
set of break points of $f_{1}$ respectively $f_{2}$. Assume that
$f_{1}$ and $f_{2}$ satisfy the conditions of Theorem \ref{DMS}. Let
$h$ be the conjugacy  between $f_{1}$ and $f_{2}$. For the proof of
Theorem \ref{DMS} we need several lemmas.

\begin{lemma}\label{lem6.1} (see \cite{DK1998}). Assume, that the
conjugating homeomorphism $h(x)$ has a positive derivative
$Dh(x_0)=\omega_0$ at some  point $x_{0}\in S^1$,
and that the following conditions hold for the points $z_i\in S^1,\,
i=1,..,4$, with $z_1\prec z_2\prec z_3\prec z_4$, and some constant
$R_{1}>1:$
\begin{itemize}
\item[(a)] the intervals $[z_{1},z_{2}],[z_{2},z_{3}],[z_{3},z_{4}]$ are pairwise
$R_{1}$-comparable;
\item[(b)] $\max\limits_{1\leq i\leq 4}\ell([z_i,x_0])\le
R_1\ell([z_1,z_2]).$
\end{itemize}
Then for any $\varepsilon>0$ there exists
$\delta=\delta(\varepsilon)>0$ such  that
\begin{eqnarray}\label{eq18}
|Dist(z_1,z_2,z_3,z_4;h)-1|\leq C_{3}\varepsilon,
\end{eqnarray}
if $z_{i}\in (x_{0}-\delta,\ x_{0}+\delta)$ for all $i=1,2,3,4$,
where the constant $C_{3}>0$ depends only on $R_{1},$ $\omega_0$ and
not on $\varepsilon.$
\end{lemma}

We define the following functions on the domain $\{(x,y):
x>0,\,0\leq y\leq 1\}$:
$$
G_{f_{1}}^{(i)}(x,y)=\frac{[\sigma_{f_{1}}^{(i)}+(1-\sigma_{f_{1}}^{(i)})y](1+x)}
{\sigma_{f_{1}}^{(i)}+(1-\sigma_{f_{1}}^{(i)})y+x},\;
1\leq i\leq  m_{1},
$$
$$
G_{f_{2}}^{(j)}(x,y)=\frac{[\sigma_{f_{2}}^{(j)}+(1-\sigma_{f_{2}}^{(j)})y](1+x)}
{\sigma_{f_{2}}^{(j)}+(1-\sigma_{f_{2}}^{(j)})y+x},\;
1\leq j\leq m_{2},
$$
where the $\sigma_{f_{1}}^{(i)}$ and $\sigma_{f_{2}}^{(j)}$ are the
jumps of $f_{1}$ and $f_{2}$  at the points $b_{1}^{(i)}$ and
$b_{2}^{(j)}$ respectively.

Let $\widehat{B}_1\subset B(f_{1})$ and $\widehat{B}_2\subset
B(f_{2}).$ Denote

$$\sigma_{f_{1}}(\widehat{B}_1):=\prod\limits_{b_{1}^{(i)}\in
\widehat{B}_1}\sigma_{f_{1}}(b_{1}^{(i)}),\,\,\sigma_{f_{2}}(\widehat{B}_2):=\prod\limits_{b_{2}^{(j)}\in
\widehat{B}_2}\sigma_{f_{2}}(b_{2}^{(j)}),$$
$$\widehat{\Lambda}_{1,2}:=\min\{\sigma_{f_{1}}(\widehat{B}_1),\,
\sigma_{f_{2}}(\widehat{B}_2),\,|\sigma_{f_{1}}(\widehat{B}_1)-\sigma_{f_{2}}(\widehat{B}_2)|\}.$$

\begin{lemma}\label{lem6.2} Let
$\widehat{B}_1=\{b_{1}^{(i_1)},b_{1}^{(i_2)},...,b_{1}^{(i_{p_1})}\}$
and
$\widehat{B}_2=\{b_{2}^{(j_1)},b_{2}^{(j_2)},...,b_{2}^{(j_{p_2})}\}$
be  arbitrary "not coinciding"  subsets of the  break points of
$f_{1}$ and $f_{2}$. Then there exist constants  $\Omega_{0}>1$ and
$\vartheta_0\in (0,1)$ such that for arbitrary
$x_{f_{1}}^{(s)},\,x_{f_{2}}^{(t)}\geq
\Omega_0,\,\,y_{f_{1}}^{(s)},\,y_{f_{2}}^{(t)}\in
[0,\vartheta_0),\,\, 1 \leq s\leq p_1,\,\,1 \leq t\leq p_2$ the
following  inequality holds:
\begin{equation}\label{eq19}
\mid\prod\limits_{s=1}^{p_1}G_{f_{1}}^{(i_s)}(x_{f_{1}}^{(s)},
y_{f_{1}}^{(s)})-\sigma_{f_{1}}(\widehat{B}_{1})\mid\leq
\frac{\widehat{\Lambda}_{1,2}}{8}
\end{equation}
\begin{equation}\label{eq20}
\mid\prod\limits_{t=1}^{p_2}G_{f_{2}}^{(j_t)}(x_{f_{2}}^{(t)},
y_{f_{2}}^{(t)})-\sigma_{f_{2}}(\widehat{B}_{2})\mid\leq
\frac{\widehat{\Lambda}_{1,2}}{8}.
\end{equation}
where $\Omega_0$ and $\vartheta_0$  only depend on
$\sigma_{f_{1}}^{(i_s)},\,1\leq s\leq p_{1}$ and
$\sigma_{f_{2}}^{(j_t)},\,1\leq t\leq p_{2}$
\end{lemma}

\textbf{Proof of Lemma \ref{lem6.2}.} Assume
$\widehat{B}_1=\{b_{1}^{(i_{s})},\,1\leq s\leq p_{1}\}$ and
$\widehat{B}_2=\{b_{2}^{(j_{t})},\,1\leq t\leq p_{2}\}$ are "not
coinciding" subsets of the break points of $f_{1}$ and $f_{2}$,
respectively. We rewrite
$\prod\limits_{s=1}^{p_1}G_{f_{1}}^{(i_s)}(x_{f_{1}}^{(s)},y_{f_{1}}^{(s)})$
in the form
$$
\prod\limits_{s=1}^{p_1}G_{f_{1}}^{(i_{s})}(x_{f_{1}}^{(s)},
y_{f_{1}}^{(s)})=\prod\limits_{s=1}^{p_1}\frac{[\sigma_{f_{1}}^{(i_{s})}+
(1-\sigma_{f_{1}}^{(i_{s})})y_{f_{1}}^{(s)}](1+x_{f_{1}}^{(s)})}
{\sigma_{f_{1}}^{(i_{s})}+(1-\sigma_{f_{1}}^{(i_{s})})y_{f_{1}}^{(s)}+
x_{f_{1}}^{(s)}}=
$$
$$\prod\limits_{s=1}^{p_1}[\sigma_{f_{1}}^{(i_{s})}+
(1-\sigma_{f_{1}}^{(i_{s})})y_{f_{1}}^{(s)}]\times\prod\limits_{s=1}^{p_1}\frac{1+x_{f_{1}}^{(s)}}
{\sigma_{f_{1}}^{(i_{s})}+
(1-\sigma_{f_{1}}^{(i_{s})})y_{f_{1}}^{(s)}+x_{f_{1}}^{(s)}}\equiv
$$
\begin{equation}\label{equa7.1}
\equiv\Phi_{f_{1}}^{(1)}(y_{f_{1}}^{(1)},...,y_{f_{1}}^{({p_1})})\times
\Phi_{f_{1}}^{(2)}(y_{f_{1}}^{(1)},...,y_{f_{1}}^{({p_1})},
x_{f_{1}}^{(1)},...,x_{f_{1}}^{({p_1})}).
\end{equation}
Obviously
$$
\lim\limits_{y_{f_{1}}^{(s)}\rightarrow 0,\,
s=\overline{1,p_1}}\Phi_{f_{1}}^{(1)}(y_{f_{1}}^{(1)},...,y_{f_{1}}^{({p_1})})=
\sigma_{f_1}(\widehat{B}_1),
$$
$$
\lim\limits_{x_{f_{1}}^{(s)}\rightarrow
\infty,\,s=\overline{1,p_1}}\Phi_{f_{1}}^{(2)}(y_{f_{1}}^{(1)},...,y_{f_{1}}^{({p_1})},
x_{f_{1}}^{(1)},...,x_{f_{1}}^{({p_1})})=1.
$$
When the variables $y_{f_{1}}^{(s)},\,\,1\leq s\leq p_1$, are
uniformly close to zero, the function
$\Phi_{f_{1}}^{(1)}=\Phi_{f_{1}}^{(1)}(y_{f_{1}}^{(1)},...,y_{f_{1}}^{({p_1})})$
hence is close to $\sigma_{f_1}(\widehat{B}_1)$, while the function
$\Phi_{f_{1}}^{(2)}=\Phi_{f_{1}}^{(2)}
(y_{f_{1}}^{(1)},...,y_{f_{1}}^{({p_1})},x_{f_{1}}^{(1)},\\...,x_{f_{1}}^{({p_1})})$
is close 1 for large values of
$x_{f_{1}}^{(s)},\,\,1\leq s\leq p_1$. Taking these remarks into
account  and using the explicit form of the functions
$\Phi_{f_{1}}^{(1)}$ and $\Phi_{f_{1}}^{(2)}$ we can now estimate
$\mid\Phi_{f_{1}}^{(1)}\cdot
\Phi_{f_{1}}^{(2)}-\sigma_{f_1}(\widehat{B}_1)\mid$. To estimate
$\Phi_{f_{1}}^{(1)}$, suppose that $0\leq
y_{f_{1}}^{(s)}\leq \vartheta_{0}^{(1)}<1$, where we shall choose
the constant $\vartheta_{0}^{(1)}$ later. It is easy to see that
$$
|\Phi_{f_{1}}^{(1)}-\sigma_{f_1}(\widehat{B}_1)|=\sigma_{f_1}
(\widehat{B}_1)|\prod\limits_{s=1}^{p_1}(1+\frac{(1-\sigma_{f_{1}}^{(i_{s})})}
{\sigma_{f_{1}}^{(i_{s})}}y_{f_{1}}^{(s)})-1|\leq
C_{4}\vartheta_{0}^{(1)},
$$
where the constant $C_{4}>0$ depends only on the
$\sigma_{f_{1}}^{(i_{s})}, 1\leq s\leq p_1$.

We set
$\vartheta_{0}^{(1)}=\min\{\frac{\widehat{\Lambda}_{1,2}}{16 C_{4}},1\}$.
Then
\begin{equation}\label{equa7.2}
\mid\Phi_{f_1}^{(1)}-\sigma_{f_1}(\widehat{B}_1)\mid<\frac{\widehat{\Lambda}_{1,2}}{16},
\end{equation}
for all $0\leq y_{f_{1}}^{(s)}\leq
\vartheta_{0}^{(1)},\,1\leq s\leq p_1$.

We next estimate $|\Phi_{f_{1}}^{(2)}-1|$ for large values of
$x_{f_{1}}^{(s)}, 1\leq s\leq p_1$. Using the explicit form of
the function $\Phi_{f_{1}}^{(2)}$, we see that the inequality
\begin{equation}\label{equa7.3}
|\Phi_{f_{1}}^{(2)}-1|<R_{2}\sum\limits_{s=1}^{p_{1}}\frac{1}{x_{f_{1}}^{(s)}},
\end{equation}
holds for all $y_{f_{1}}^{(s)},\,0\leq y_{f_{1}}^{(s)}\leq 1$, and
$x_{f_{1}}^{(s)}>0, 1\leq s\leq p_1$, where the constant $R_2>0$
depends only on
$\sigma_{f_{1}}^{(i_{s})}, 1\leq s\leq p_1$.

Suppose now, that $x_{f_{1}}^{(s)}\geq \Omega_{0}^{(1)}, 1\leq s\leq
p_1$. By (\ref{equa7.3}) we have
$|\Phi_{f_1}^{(2)}-1|<R_{3}p_{1}\frac{1}{\Omega_{0}^{(1)}}$. It
follows from this together with (\ref{equa7.2}) that
$|\Phi_{f_{1}}^{(1)}\cdot\Phi_{f_{1}}^{(2)}-\sigma_{f_1}(\widehat{B}_1)|\leq
|\Phi_{f_{1}}^{(1)}-\sigma_{f_1}(\widehat{B}_1)|+|\Phi_{f_{1}}^{(1)}|\cdot|\Phi_{f_{1}}^{(2)}-1|\leq
\frac{\widehat{\Lambda}_{1,2}}{16}+(\sigma_{f_1}(\widehat{B}_1)+
\frac{\widehat{\Lambda}_{1,2}}{16})R_{3}p_1\frac{1}{\Omega_{0}^{(1)}}$.
We choose $\Omega_{0}^{(1)}$ in such a way that the relation
$\frac{\widehat{\Lambda}_{1,2}}{16}+(\sigma_{f_1}(\widehat{B}_1)+
\frac{\widehat{\Lambda}_{1,2}}{16})R_{3}p_1\frac{1}{\Omega_{0}^{(1)}}=\frac{\widehat{\Lambda}_{1,2}}{8}$
holds, whence, $\Omega_{0}^{(1)}=\frac{16\sigma_{f_1}(\widehat{B}_1)+\widehat{\Lambda}_{1,2}}{\widehat{\Lambda}_{1,2}}R_{3}p_1$.
As a result we have
$|\Phi_{f_1}^{(1)}\cdot\Phi_{f_1}^{(2)}-\sigma_{f_1}(\widehat{B}_1)|\leq
\frac{\widehat{\Lambda}_{1,2}}{8}$ for $0\leq y_{f_{1}}^{(s)}\leq
\vartheta_{0}^{(1)}$ and $x_{f_{1}}^{(s)}\geq \Omega_{0}^{(1)},
1\leq s\leq p_1$. It follows from this  that the assertion
(\ref{eq19}) of the lemma holds. Analogously it can be shown that with
\begin{equation}\label{equa7.4}
\vartheta_{0}^{(2)}:=\min\{\frac{\widehat{\Lambda}_{1,2}}{16C_{5}},1\},\,\,\Omega_{0}^{(2)}:=
\frac{16\sigma_{f_{2}}(\widehat{B}_2)+\widehat{\Lambda}_{1,2}}{\widehat{\Lambda}_{1,2}}R_{4}p_2,
\end{equation}
and $0\leq y_{f_{2}}^{(t)}\leq
\vartheta_{0}^{(2)}$ and $x_{f_{2}}^{(t)}\geq \Omega_{0}^{(2)},
1\leq t\leq p_2,$  the assertion (\ref{eq20}) of Lemma
\ref{lem6.2} holds. In (\ref{equa7.4}) the constants $C_5>0$ and
$R_4>0$ depend on the $\sigma_{f_{2}}^{ (j_t)}, 1\leq t\leq p_2$. If
we finally set
$\vartheta_{0}:=\min\{\vartheta_{0}^{(1)},\vartheta_{0}^{(2)}\}$
and $\Omega_{0}:=\max\{\Omega_{0}^{(1)},\Omega_{0}^{(2)}\}$  Lemma
\ref{lem6.2} holds for $x_{f_{1}}^{(s)},\,x_{f_{2}}^{(t)}\geq
\Omega_0$ and $y_{f_{1}}^{(s)},\,y_{f_{2}}^{(t)}\in
[0,\vartheta_0),\,\,1\leq s\leq p_1,\, 1\leq t\leq p_2$.

Define 
\begin{equation}\label{equaMZ1}
\overline{\Omega}_0:=\max\Omega_0(\sigma_{f_{1}}^{(i_1)},...,\sigma_{f_{1}}^{(i_{p_1})},
\sigma_{f_{2}}^{(j_1)},...,\sigma_{f_{2}}^{(j_{p_2})})
\end{equation}
\begin{equation}\label{equaMZ2}
\overline{\vartheta}_0:=\min\vartheta_0(\sigma_{f_{1}}^{(i_1)},...,\sigma_{f_{1}}^{(i_{p_1})},
\sigma_{f_{1}}^{(j_1)},...,\sigma_{f_{2}}^{(j_{p_2})})
\end{equation}
where the minimum and maximum are taken over all "not coinciding"
subsets $\widehat{B}_1$ and $\widehat{B}_2$ of the break points of
$f_1$ and $f_2$ and $v_1,v_2>0$ are the total variations of $\ln
Df_{1}$ and $\ln Df_{2}$ over $S^{1}$ respectively. Next we define
the following constants $M_0$ and $\zeta_0$:
$$
M_0:=\overline{\Omega}_0 e^{\max\{v_1,v_2\}},\,\,
\zeta_0:=\overline{\vartheta}_0 e^{-\min\{v_1,v_2\}}.
$$
Let $K_0=K_0(f_1,f_2,M_0,\zeta_0)>M_0>1$ be the constant $K$ as
defined in Theorem \ref{CO}.
\begin{lemma}\label{lem6.3} Suppose the circle homeomorphisms $f_{1}$ and $f_{2}$ satisfy
the conditions of Theorem \ref{DMS} and
$Dh(x_{0})=\omega_{0}>0$ at some point $x_{0}\in
S^{1}.$ Let be $\delta>0$ and assume the triples intervals
$([z_s,z_{s+1}],\,\,s=1,2,3)$ and
$([h(z_s),h(z_{s+1})],\,\,s=1,2,3)$\,\, cover the  "not coinciding"
subsets $\widehat{B}_1$ and $\widehat{B}_2$ of break points of
$f_{1}$ and $f_{2}$ "$(K_0,M_0,\delta,\zeta_{0};x_{0},h(x_0))$-
regularly" for some $r_{n}\in \{q_{n-1},q_{n},q_{n}+q_{n-1}\}$. Then
for sufficiently large $n$ the following inequality holds:
\begin{equation}\label{equadist}
\mid\frac{Dist(z_1,z_2,z_3,z_4;f_{1}^{r_n})}{Dist(h(z_1),h(z_2),h(z_3),h(z_4);f_{2}^{r_n})}-1\mid\geq
R_5>0,
\end{equation}
where the constant $R_5$  depends only on $f_1$ and $f_2$.
\end{lemma}

\textbf{Proof of Lemma \ref{lem6.3}.}  Suppose, the triples of
intervals $([z_s,z_{s+1}],\,\,s=1,2,3)$ and
$([h(z_s),h(z_{s+1})],\,\,s=1,2,3)$
cover the "not
coinciding" subsets $\widehat{B}_1$ and $\widehat{B}_2$ of break
points of $f_1$ and $f_2$ "$(K_0,M_0,\delta,\zeta_0;x_0,h(x_0))$-regularly" for some $r_n\in
\{q_{n-1},q_{n},q_{n-1}+q_{n-1}\}$. To be definite, assume
$\widehat{B}_1=\{b_{1}^{(i_1)},b_{1}^{(i_2)},...,b_{1}^{(i_{p_1})}\}\neq \O$,\,
$\widehat{B}_2=\{b_{2}^{(j_1)},b_{2}^{(j_2)},...,b_{2}^{(j_{p_2})}\}\neq \O$.
 By Definition \ref{defCO1}, the triples of  intervals
$([z_s,z_{s+1}],\,\,s=1,2,3)$ and
$([h(z_s),h(z_{s+1})],\,\,s=1,2,3)$
cover the  subsets
$\widehat{B}_1$ respectively $\widehat{B}_2$ "$(K_0,M_0,\delta,\zeta_0;x_0)$-regularly" respectively
"$(K_0,M_0,\delta,\zeta_0;h(x_0))$-regularly".

By Definition \ref{defCO} we have $\overline{b}_{1}^{(i_{s})}\in
[z_1,z_2],\,\,1 \leq s\leq p_1$ and $\overline{b}_{2}^{(j_{t})}\in
[h(z_1),h(z_2)],\,\,1\leq t\leq p_2$. Notice that the
intervals
$f_{1}^{l_{1}^{(i_{s})}}([z_{1},z_{2}]),\,1\leq s\leq p_1$
cover the break points $b_{1}^{(i_{s})},\,1\leq s\leq p_1$.
 Similarly, the intervals
$f_{2}^{l_{2}^{(j_{t})}}([h(z_{1}),h(z_{2})]),\,1\leq t\leq p_2$
cover the break points $b_{2}^{(j_{t})},\,1\leq t
\leq p_2.$  Next we want to compare the distortion $Dist(z_1, z_2,
z_3, z_4; f_{1}^{r_n})$ and \\ $Dist(h(z_1),h(z_2),h(z_3),h(z_4);
f_{2}^{r_n}).$ We estimate only the first distortion, the second one
can be estimated analogously. Obviously
$$
Dist(z_1, z_2, z_3, z_4; f_{1}^{r_n})=\frac{Cr(f_{1}^{r_{n}}(z_1),
f_{1}^{r_{n}}(z_2), f_{1}^{r_{n}}(z_3), f_{1}^{r_{n}}(z_4))}{Cr(z_1,
z_2, z_3, z_4)}=$$
$$
=\frac{Cr(f_{1}(z_1), f_{1}(z_2), f_{1}(z_3), f_{1}(z_4))}{Cr(z_1,
z_2, z_3, z_4)}\times\frac{Cr(f_{1}^{2}(z_1), f_{1}^{2}(z_2),
f_{1}^{2}(z_3), f_{1}^{2}(z_4))}{Cr(f(z_1), f(z_2), f(z_3),
f(z_4))}\times...$$
$$
\times\frac{Cr(f_{1}^{r_{n}}(z_1), f_{1}^{r_{n}}(z_2),
f_{1}^{r_{n}}(z_3), f_{1}^{r_{n}}(z_4))}{Cr(f_{1}^{r_{n}-1}(z_1),
f_{1}^{r_{n}-1}(z_2), f_{1}^{r_{n}-1}(z_3),
f_{1}^{r_{n}-1}(z_4))}=\prod\limits_{i=0}^{r_{n}-1} Dist
(f_{1}^{i}(z_1), f_{1}^{i}(z_2), f_{1}^{i}(z_3),
f_{1}^{i}(z_4);f_{1})
.$$

We rewrite $Dist(z_1, z_2, z_3, z_4; f_{1}^{r_n})$ in the form
$$
Dist (z_1, z_2, z_3, z_4; f_{1}^{r_n})=\prod\limits_{s=1}^{p_1} Dist
(f_{1}^{l_{1}^{(i_s)}}(z_1), f_{1}^{l_{1}^{(i_s)}}(z_2),
f_{1}^{l_{1}^{(i_s)}}(z_3), f_{1}^{l_{1}^{(i_s)}}(z_4);f_{1})\times
$$
\begin{equation}\label{equa24}
\times \prod\limits_{i=0,i\not=l_{1}^{(i_s)},
s=\overline{1,p_1}}^{r_n-1}Dist (f_{1}^i(z_1), f_{1}^i(z_2),
f_{1}^i(z_3), f_{1}^i(z_4);f_{1}).
\end{equation}
To estimate the first factor in (\ref{equa24}) we use Lemma
\ref{DAKX2} and the definition of the functions
$G_{f_{1}}^{(i)}(x,y)$ to get
$$
Dist (f_{1}^{l_{1}^{(i_s)}}(z_1), f_{1}^{l_{1}^{(i_s)}}(z_2),
f_{1}^{l_{1}^{(i_s)}}(z_3), f_{1}^{l_{1}^{(i_s)}}(z_4); f_{1})=
$$
$$
=\frac{[\sigma_{f_{1}}^{(i_s)}+(1-\sigma_{f_{1}}^{(i_s)})z_{f_{1}}^{(i_s)}(l_{1}^{(i_s)})
] (1+\xi_{f_{1}}(l_{1}^{(i_s)}))}{\sigma_{f_{1}}^{(i_s)}+
(1-\sigma_{f_{1}}^{(i_s)})z_{f_{1}}^{(i_s)}(l_{1}^{(i_s)})+\xi_{f_{1}}(l_{1}^{(i_s)})}+\chi_{s}^{(1)}=
$$
\begin{equation}\label{equa23}
=G_{f_{1}}^{(i_s)}(\xi_{f_{1}}(l_{1}^{(i_s)}),
z_{f_{1}}^{(i_s)}(l_{1}^{(i_s)}))+\chi_{s}^{(1)},\quad
1 \leq s\leq p_1,
\end{equation}
where $|\chi_{s}^{(1)}|\leq C_2
\ell([f_{1}^{l_{1}^{(i_{s})}}(z_{1}),f_{1}^{l_{1}^{(i_{s})}}(z_{4})]),\,1\leq s\leq p_1$.
By construction the interval $[z_1,z_4]$ (see
(\ref{Z1}),(\ref{Z2}),(\ref{Z3})) is $r_n$-small and therefore the
intervals $f_{1}^{i}([z_1,z_4]),\,0<i<r_n$ are pairwise disjoint.
Hence, using Corollary \ref{LDB} we obtain
\begin{equation}\label{ellf}
\ell(f_{1}^{i}([z_{1},z_{4}]))\leq
const\lambda_{1}^{n},\,i=\overline{0,r_{n}-1},
\end{equation}
where $\lambda_{1}=(1+e^{-v_{1}})^{-\frac{1}{2}}<1.$

Fix now some $\varepsilon>0$. There exists $N=N(\varepsilon)>1$ such
that \begin{equation}\label{equa25}
|\chi_{s}^{(1)}|<C_6\varepsilon,\,1\leq s\leq p_1
\end{equation}
 holds for $n>N$, where the constant $C_{6}>0$ depends only on $f_{1}$.

Suppose, $\xi_{f_{1}}(0)$ and $z_{f_{1}}^{(i_{s})}(0),\,1\leq s\leq
p_1$ satisfy the following conditions: $\xi_{f_{1}}(0)>M_0$ and
$z_{f_{1}}^{(i_{s})}(0)<\zeta_0$ for $1\leq s\leq p_1$. Then, using
relations (\ref{evf}) we get
$\xi_{f_{1}}(l_{1}^{(i_s)})>\overline{\Omega}_{0}$ and
$z_{f_{1}}^{(i_s)}(l_{1}^{(i_s)})<\overline{\vartheta}_{0},\,\,1\leq
s\leq p_1$, where $\overline{\Omega}_{0}$ and
$\overline{\vartheta}_{0}$ defined in (\ref{equaMZ1}) and
(\ref{equaMZ2}). Since $\overline{\vartheta}_{0}$
is the minimum  of $\vartheta_{0}$, the assertion of Lemma \ref{lem6.2}
is true for $\overline{\vartheta}_{0}$ also. It follows from the
assertion of Lemma \ref{lem6.2} that
\begin{equation}\label{equa26}
|\prod\limits_{s=1}^{p_1}G_{f_{1}}^{(i_{s})}(\xi_{f_{1}}(l_{1}^{(i_{s})}),
z_{f_{1}}^{(i_{s})}(l_{1}^{(i_{s})}))-\sigma_{f_{1}}(\widehat{B}_1)|\leq
\frac{\widehat{\Lambda}_{1,2}}{8}.
\end{equation}
By combining (\ref{equa23})-(\ref{equa26}) we obtain
\begin{equation}\label{equa27}
|\prod\limits_{s=1}^{p_1}[G_{f_{1}}^{(i_{s})}(\xi_{f_{1}}(l_{1}^{(i_{s})}),
z_{f_{1}}^{(i_{s})}(l_{1}^{(i_{s})}))+\chi_{s}^{(1)}]-\sigma_{f_1}(\widehat{B}_1)|\leq
\frac{\widehat{\Lambda}_{1,2}}{6}
\end{equation}
for sufficiently small $\varepsilon>0$.

Next we estimate the second factor in (\ref{equa24}).  Applying
Lemma \ref{DAKX1} we obtain
$$
\prod\limits_{i=0, i\not=l_{1}^{(i_s)},
s=\overline{1,p_1}}^{r_n-1}Dist (f_{1}^i(z_1),f_{1}^i(z_2),
f_{1}^i(z_3), f_{1}^i(z_4);f_{1}) =
$$
\begin{equation}\label{equa28}
=\exp\{\sum\limits_{i=0,i\not=l_{1}^{(i_s)},
s=\overline{1,p_1}}^{r_{n}-1}\log
(1+O((\ell([f_{1}^{i}(z_{1}),f_{1}^{i}(z_{4})]))^{1+\alpha}))\}
\end{equation}
 Using the bound (\ref{ellf}) and  that the  interval $[z_{1},z_{4}]$ is
$r_{n}$-small, we obtain  from (\ref{equa28})
$$
|\prod\limits_{i=0, i\not=l_{1}^{(i_s)},
s=\overline{1,p_1}}^{r_n-1}Dist (f_{1}^i(z_1),
 f_{1}^i(z_2),
f_{1}^i(z_3), f_{1}^i(z_4); f_{1})-1|\leq $$
\begin{equation}\label{equa29}
\leq const \lambda_{1}^{n\alpha}
\sum\limits_{i=0,i\not=l_{1}^{(i_s)},
s=\overline{1,p_1}}^{r_{n}-1}(\ell([f_{1}^{i}(z_{1}),f_{1}^{i}(z_{4})])\leq
const \lambda_{1}^{n\alpha}.
\end{equation}
The relations (\ref{equa27}) and (\ref{equa29}) imply that for
sufficiently large $n$
\begin{equation}\label{equa30}
|Dist(z_1, z_2, z_3, z_4;
f_{1}^{r_n})-\sigma_{f_{1}}(\widehat{B}_1)|<\frac{\widehat{\Lambda}_{1,2}}{4}
\end{equation}
holds.

The same way  it can be  shown that
 for the triple of  intervals
$([h(z_s),h(z_{s+1})],\,\,s=1,2,3)$ covering the set $\widehat{B}_2$ "$(K_0,M_0,\delta,\zeta_0;h(x_0))$-regularly"  the following inequality
\begin{equation}\label{equa31}
|Dist(h(z_1),h(z_2),h(z_3),h(z_4);
f_{2}^{r_n})-\sigma_{f_{2}}(\widehat{B}_2)|<\frac{\widehat{\Lambda}_{1,2}}{4}
\end{equation}
holds for
sufficiently large $n$ . The inequalities (\ref{equa30}) and (\ref{equa31}) show that
\begin{equation}\label{equa32}
\frac{Dist(z_1, z_2, z_3, z_4; f_{1}^{r_n})}{Dist(h(z_1),h(z_2),
h(z_3),h(z_4); f_{2}^{r_n})}-1\geq
\frac{4(\sigma_{f_{1}}(\widehat{B}_1)-\sigma_{f_{2}}(\widehat{B}_2))-
2\widehat{\Lambda}_{1,2}}{4\sigma_{f_{2}}(\widehat{B}_2)
+\widehat{\Lambda}_{1,2}}>0,
\end{equation}
 if   $\sigma_{f_{1}}(\widehat{B}_1)>\sigma_{f_{2}}(\widehat{B}_2)$, and
\begin{equation}\label{equa33}
\frac{Dist(z_1, z_2, z_3, z_4;
f_{1}^{r_n})}{Dist(h(z_1),h(z_2),h(z_3),h(z_4); f_{2}^{r_n})}-1\leq
\frac{4(\sigma_{f_{1}}(\widehat{B}_1)-\sigma_{f_{2}}(\widehat{B}_2))+
2\widehat{\Lambda}_{1,2}}{4\sigma_{f_{2}}(\widehat{B}_2)
-\widehat{\Lambda}_{1,2}}<0,
\end{equation}
if
$\sigma_{f_{1}}(\widehat{B}_1)<\sigma_{f_{2}}(\widehat{B}_2)$.
 If we set
\begin{equation}\label{equa34}
R_{5}:=\min\{\frac{|4(\sigma_{f_{1}}(\widehat{B}_1)-\sigma_{f_{2}}(\widehat{B}_2))-
2\widehat{\Lambda}_{1,2}|}{4\sigma_{f_{2}}(\widehat{B}_2)+\widehat{\Lambda}_{1,2}},\,
\frac{|4(\sigma_{f_{1}}(\widehat{B}_1)-\sigma_{f_{2}}(\widehat{B}_2))+
2\widehat{\Lambda}_{1,2}|}{4\sigma_{f_{2}}(\widehat{B}_2)-\widehat{\Lambda}_{1,2}}
\},
\end{equation}
where the minimum is taken over all "not coinciding" subsets
$\widehat{B}_1$ and $\widehat{B}_2$ of break points $f_1$ and $f_2$, then
It follows from (\ref{equa32})-(\ref{equa34}) that the assertion of
the lemma holds.

\textbf{Proof of Theorem \ref{DMS}.} Let $f_{1}$ and $f_{2}$ be
circle homeomorphisms satisfying the conditions of Theorem
\ref{DMS}. The lift $H(x)$ of the conjugating map $h(x)$ is a
continuous and monotone increasing function on $R^1$. Hence $H(x)$
has a finite derivative $DH(x)$ for almost all $x$
with respect to Lebesgue measure. We claim that
$Dh(x)=0$ at all the points $x$ where the finite
derivative exists. Suppose $Dh(x_0)>0$ for some
point $x_0\in S^1$. We fix $\varepsilon>0$. Let
$K_0=K_0(f_1,f_2,M_0,\zeta_0)>M_0>1$ be the constant defined in the
assertion of Theorem \ref{CO}. By the assertion  of Theorem
\ref{CO}, for sufficiently large $n$ there exist  "not coinciding"
subsets $\widehat{B}_1$ and $\widehat{B}_2$ of break points of $f_1$
and $f_2$, points $z_{i}\in S^{1},\,1\leq i\leq 4,$ with $z_{1}\prec
z_{2}\prec z_{3}\prec z_{4}$ and a number $r_{n}\in
\{q_{n-1},q_{n},q_{n}+q_{n-1}\}$ such that  the triples of intervals
$[z_s,z_{s+1}],\,\,s=1,2,3,$ and $[h(z_s),h(z_{s+1})],\,\,s=1,2,3,$
 cover the
points of $\widehat{B}_1$ and $\widehat{B}_2$ "$(K_{0},M_{0},\delta,\zeta_{0};x_{0},h(x_0))$-regularly" for $r_{n}$. By
Definition \ref{defCO} of a regularly covering each of the systems of   intervals
$([z_s,z_{s+1}],\,\,s=1,2,3)$ and
$([f_{1}^{r_n}(z_s),f_{1}^{r_n}(z_{s+1})],\,\,s=1,2,3)$  satisfies
the conditions of Lemma \ref{lem6.1} with constant $R_1=K_0.$

Using the assertion of Lemma \ref{lem6.1} we obtain
\begin{equation}\label{equa21}
|Dist(z_1,z_2,z_3,z_4;h)-1|\leq C_{3}\varepsilon,
\end{equation}
\begin{equation}\label{equa22}
|Dist(f_{1}^{r_n}(z_1),f_{1}^{r_n}(z_2),f_{1}^{r_n}(z_3),f_{1}^{r_n}(z_4);h)-1|\leq
C_{3}\varepsilon.
\end{equation}
Hence
\begin{equation}\label{equa21}
|\frac{Dist(z_1,z_2,z_3,z_4;h)}{Dist(f_{1}^{r_n}(z_1),f_{1}^{r_n}(z_2),
f_{1}^{r_n}(z_3),f_{1}^{r_n}(z_4);h)}-1|\leq
C_7\varepsilon,
\end{equation}
where the constant $C_7>0$ does not depend on
$\varepsilon$ and $n$.

Since $h$ is conjugating $f_{1}$ and $f_{2}$ we can readily see that
$$
Cr(h(f_{1}^{r_n}(z_1)),h(f_{1}^{r_n}(z_2)),h(f_{1}^{r_n}(z_3)),h(f_{1}^{r_n}(z_4)))=
$$
$$
=Cr(f_{2}^{r_n}(h(z_1)),f_{2}^{r_n}(h(z_2)),f_{2}^{r_n}(h(z_3)),f_{2}^{r_n}(h(z_4))).
$$
Hence we obtain
$$
\frac{Dist(f_{1}^{r_n}(z_1),f_{1}^{r_n}(z_2),f_{1}^{r_n}(z_3),
f_{1}^{r_n}(z_4);h)}{Dist(z_1,z_2,z_3,z_4;h)}=
$$
$$
=\frac{Cr(h(f_{1}^{r_n}(z_1)),h(f_{1}^{r_n}(z_2)),h(f_{1}^{r_n}(z_3)),
h(f_{1}^{r_n}(z_4)))}{Cr(f_{1}^{r_n}
(z_1),f_{1}^{r_n}(z_2),f_{1}^{r_n}(z_3),f_{1}^{r_n}(z_4))}\times
$$
$$
\times \frac{Cr(z_1,z_2,z_3,z_4)}{Cr(h(z_1),h(z_2),h(z_3),h(z_4))}=
\frac{Cr(f_{2}^{r_n}(h(z_1)),f_{2}^{r_n}(h(z_2)),f_{2}^{r_n}(h(z_3)),
f_{2}^{r_n}(h(z_4)))}{Cr(h(z_1),h(z_2),h(z_3),h(z_4))}:
$$
$$
:\frac{Cr(f_{1}^{r_n}(z_1),f_{1}^{r_n}(z_2),f_{1}^{r_n}(z_3),
f_{1}^{r_n}(z_4))}{Cr(z_1,z_2,z_3,z_4)}=
\frac{Dist(h(z_1),h(z_2),h(z_3),h(z_4);f_{2}^{r_n})}
{Dist(z_1,z_2,z_3,z_4;f_{1}^{r_n})}.
$$
This, together with (\ref{equa21}) obviously implies that
$$
\mid\frac{Dist(z_1,z_2,z_3,z_4;f_{1}^{r_n})}
{Dist(h(z_1),h(z_2),h(z_3),h(z_4);f_{2}^{r_n})}-1\mid\leq
C_8\varepsilon.
$$
where the constant $C_8>0$ does not depend on
$\varepsilon$ and $n$. But this contradicts equation
(\ref{equadist}). Theorem \ref{DMS} is
 therefore completely proved.

\end{document}